\documentclass[final,5p,times,twocolumn]{elsarticle}

\journal{Systems \& Control Letters}

\usepackage[numbers]{natbib}
\usepackage[utf8]{inputenc} 
\usepackage[english]{babel} 
\usepackage{amsmath,amssymb,mathtools, comment} 
\usepackage{graphicx} 
\graphicspath{{images/}} 
\usepackage{float}

\usepackage{Macros}

\newproof{proof}{Proof}
\newproof{pot}{Proof of Theorem \ref{thm}}

\newtheorem{lemma}{Lemma}[section]

\newtheorem{proposition}[lemma]{Proposition}
\newtheorem{theorem}[lemma]{Theorem}
\newtheorem{assumption}[lemma]{Assumption}

\newtheorem{definition}[lemma]{Definition}
\newtheorem{example}[lemma]{Example}

\newtheorem{remark}[lemma]{Remark}

\renewcommand\theenumi{(\roman{enumi})}

\allowdisplaybreaks[4]  

    \usepackage{bibentry}
    \nobibliography*

\usepackage{tikz}
\usetikzlibrary{positioning}
\usepackage{pgfplots}
\pgfplotsset{compat=1.18}
\pgfplotsset{soldot/.style={color=uniwueblue,only marks,mark=*, mark options={scale=.3mm}}}
\pgfplotsset{holdot/.style={color=uniwueblue,fill opacity=0, only marks,mark=*, mark options={scale=.3mm}}}
\usetikzlibrary{cd}
\usetikzlibrary{shapes,fit} 

\usetikzlibrary{decorations.markings}
\tikzset{degil/.style={
		decoration={markings,
			mark= at position 0.5 with {
				\node[transform shape] (tempnode) {$\backslash$};
			}
		},
		postaction={decorate}
	}
}

\usepackage[final]{hyperref}			
\hypersetup{
    colorlinks = true,
    allcolors=blue,
}
\usepackage{mathtools}
\begin{document}
\begin{frontmatter}




\title{Converse Lyapunov theorem for input-to-state stability
of infinite-dimensional systems\tnoteref{dedication}}


\tnotetext[dedication]{To Eduardo Sontag on the occasion
of his 75th birthday with gratitude for his foundational contributions to control theory}

\author[1]{Patrick Bachmann\fnref{fn1}}
\ead{patrick.bachmann@uni-wuerzburg.de}
%
%
\affiliation[1]{organization={Institute of Mathematics, Julius-Maximilians-Universität Würzburg},
                addressline={Emil-Fischer-Str. 40},
                postcode={97074}, 
                postcodesep={}, 
                city={Würzburg},
                country={Germany}}

\author[2]{Andrii Mironchenko\corref{cor1}\fnref{fn1}}
\ead{andrii.mironchenko@uni-bayreuth.de}
\ead[URL]{http://mironchenko.com}


\affiliation[2]{organization={Department of Mathematics, University of Bayreuth},
                postcode={95440}, 
                postcodesep={}, 
                city={Bayreuth},
                country={Germany}}

\cortext[cor1]{Corresponding author}
\fntext[fn1]{P. Bachmann and A. Mironchenko have been supported by the German Research Foundation (DFG), grants Nr. 552243517 (DA 767/18-1) and 540580186 (MI 1886/3-1), respectively.}

\begin{abstract}
We prove converse Lyapunov theorems for input-to-state stability (ISS) of a broad class of infinite-dimensional systems applicable to evolution equations in Banach spaces, boundary control systems, and delay systems. 
To this end, we employ a novel method that does not require constructing auxiliary closed-loop systems. 
We start by showing that ISS of a control system implies uniform global asymptotic stability of the same system with respect to trajectory-dominated inputs. Then, under an additional regularity requirement, we prove that an ISS Lyapunov function can be constructed. 

Based on our general result, we prove that for linear systems with unbounded $\infty$-admissible operators, ISS is equivalent to the existence of a Lipschitz continuous coercive ISS Lyapunov function. 
\end{abstract}







\end{frontmatter}


\section{Introduction}

The concept of input-to-state stability (ISS), introduced by E. Sontag in the late 1980s in \cite{Son89}, unified the Lyapunov and input-output stability theories and revolutionized the constructive nonlinear control theory for  ordinary differential equations (ODEs) \cite{KoA01}.
It played a major role in robust stabilization of
nonlinear systems \cite{FrK08}, design of robust 
nonlinear observers \cite{Lib07, AnP09}, nonlinear
detectability \cite{SoW97,KSW01}, event-based control \cite{Tab07,PTN15}, stability of nonlinear networks \cite{JMW96, DRW07},
supervisory adaptive control \cite{HeM99b}, data-driven control \cite{CBD25,PPT24,LaA23}, and other areas.

The success of the ISS theory for ODEs, together with the demand for suitable tools for the robust stability analysis of partial differential equations (PDEs), has stimulated the development of the ISS theory of infinite-dimensional systems \cite{MiP20,Sch20,KaK19}. 


From the early days of the infinite-dimensional ISS theory, an abstract point of view has been adopted \cite{KaJ07,DaM13}, where control systems are defined as a flow map satisfying some basic properties. This development has led to the proof of the general direct ISS Lyapunov theorems \cite{KaJ07,DaM13}, ISS superposition theorems characterizing ISS in terms of stability and (weak) attractivity \cite{MiW18b,BDM26}, ISS theory of linear systems \cite{Sch20,JNPS18,APS26}, small-gain theory for finite and infinite networks, non-coercive Lyapunov theory \cite{MiW18b,JMP20}, etc. 
These abstract results can serve as a basis for the development of ISS theories specialized for 
parabolic PDEs, hyperbolic PDEs, delay systems, ODEs, infinite networks, and a variety of other system classes and important physical equations. 


\textbf{Lyapunov methods.} Constructing an ISS Lyapunov function is the most realistic approach to prove ISS of a given infinite-dimensional system. Lyapunov methods have been applied to nonlinear parabolic systems with distributed disturbances \cite{MaP11,DaM13}, linear time-varying systems of conservation laws \cite{PrM12}, and parabolic systems with Robin and Neumann boundary inputs \cite{ZhZ18,Sch20}. 

In contrast to this, no construction of coercive ISS Lyapunov functions has been provided for the linear 1-D heat equation with Dirichlet boundary inputs, which is exponentially stable in the absence of inputs.
This failure to apply Lyapunov methods to a toy parabolic example has  led to the development of non-Lyapunov methods for ISS analysis: admissibility theory \cite{JNPS18} and spectral analysis \cite{KaK16,LhS19} have been successfully used for linear systems; a monotonicity approach \cite{MKK19} can be applied to monotone control systems such as nonlinear parabolic systems with distributed and boundary inputs; De Giorgi iteration \cite{ZhZ19} has been applied for (local) ISS analysis of Burgers' equation.
Despite the efficiency of these methods within particular system classes, they lack the generality of the Lyapunov method.

To address the problem with the applicability of Lyapunov tools, in \cite{MiW18b} the concept of a non-coercive ISS Lyapunov function --- i.e.\ a Lyapunov function that is not necessarily coercive --- has been proposed, and in \cite{JMP20}, it was shown that the existence of a non-coercive ISS Lyapunov function implies ISS, provided some additional properties of the system are satisfied. In particular, in \cite{JMP20}, a quadratic, non-coercive ISS Lyapunov function was constructed for the 1-D heat equation with a Dirichlet boundary input. 
In \cite{Mar26}, a construction of coercive Lyapunov functions was proposed for systems that are regular, possess a non-coercive Lyapunov function, and for which there exists an exactly observable output. However, parabolic systems with Dirichlet inputs do not belong to this class. 
Another relaxation of the Lyapunov function concept, called generalized ISS Lyapunov function, has been proposed in \cite{ZhZ26}, which gives a useful and flexible tool for Lyapunov analysis of ISS for parabolic and hyperbolic PDEs with boundary disturbances, but only direct Lyapunov results for this type of Lyapunov functions are available.

\textbf{Converse ISS Lyapunov theory} was initiated by Sontag and Wang in their celebrated paper \cite{SoW95}, where it was shown that for ODEs with Lipschitz nonlinearities ISS is equivalent to the existence of a smooth ISS Lyapunov function. 
The key idea in the proof technique of \cite{SoW95} is to reduce the question of ISS of the original system to the analysis of uniform global asymptotic stability (UGAS) of an auxiliary system, obtained by closing the loop via a disturbed feedback. The latter problem can be solved using general converse Lyapunov results for UGAS developed in \cite{LSW96}. 

The approach from \cite{SoW95} has motivated numerous Lyapunov characterizations for ISS-like properties for various classes of systems with inputs \cite{Son08}. In particular, for time-delay systems, converse ISS Lyapunov theorems have been reported in  \cite{KPJ08}, and 
for evolution equations with Lipschitz nonlinearities and distributed controls, ISS Lyapunov characterizations  have been shown in \cite{MiW17c}. 

Converse Lyapunov results for boundary control systems are harder to tackle. 
Recently, \cite{MiS26} showed that linear analytic ISS systems on a Hilbert space in general do not possess coercive quadratic ISS Lyapunov functions, even if the input operator is bounded. However, as argued in \cite{MiS26}, if the input operator is $p$-admissible with $p<2$, and the semigroup is analytic and similar to a contraction, then a coercive quadratic $L^2$-ISS Lyapunov function always exists. This clarifies the difference between the heat equation with Neumann boundary conditions and with Dirichlet boundary conditions. In the first case, the admissibility properties are better, and a coercive quadratic ISS Lyapunov function exists, while in the second case one must sacrifice either coercivity or quadratic structure. 
Nevertheless, \cite{MiS26} does not rule out the existence of general non-quadratic ISS Lyapunov functions for linear systems for which quadratic ones do not exist.

\textbf{Contribution.} 
In this work, \emph{our main objective is to prove a converse ISS Lyapunov theorem for a general class of infinite-dimensional control systems, encompassing  nonlinear evolution equations with Lipschitz nonlinearities, nonlinear time-delay systems, as well as boundary control systems.}

The classical approach from \cite{SoW95,MiW17c} was formulated for systems with Lipschitz continuous right-hand sides/nonlinearities, and does not apply to general control systems for several reasons. Firstly, linear systems with unbounded input operators do not have Lipschitz continuous (or even continuous) right-hand sides. Second, abstract systems given by the flow have no right-hand sides at all. 
Thirdly, constructing auxiliary systems by closing the loop via disturbed feedback for systems given by the flow without equations defining the dynamics is hard to implement.

To overcome these obstructions, we revisit our proof strategy.
Having defined the general class of systems and families of admissible inputs for a control system in Section~\ref{sec:Abstract control systems}, in Section~\ref{sec:Discont-converseLyapunovThm}
we proceed to showing that UGAS with respect to this family of inputs is equivalent to the existence of a discontinuous Lyapunov function certifying UGAS.
The key point in this result is that we do not require any extra regularity assumptions from the flow of the system. 
The authors have not found such a result in the continuous-time setting in the literature. However,  constructions of this type are well-known for discrete-time systems, where continuity of a Lyapunov function is not essential \cite{GrK14}.
Next, assuming certain regularity of the flow map, we construct Lipschitz continuous UGAS Lyapunov functions for UGAS systems using discontinuous Lyapunov functions as building blocks, using known methods \cite{KaJ11a,Hen81}.

Having developed converse Lyapunov theory for UGAS systems at a level of generality needed for our purposes, we proceed in Section \ref{sec:Converse ISS Lyapunov theorem} to the development of converse Lyapunov results for ISS. 
\emph{We prove that ISS implies UGAS with respect to trajectory-dominated inputs}, i.e., inputs satisfying 
\[
\|u(t)\|_U\leq \eta(\|\phi(t,x,u)\|_X),
\]
for a certain $\eta \in\Kinf$.
One can view $\eta$ as a robustness margin of a control system, showing how big an input can be with respect to the state of the system, such that the system remains UGAS for all such inputs. 
This transformation is important as it is hard to construct Lyapunov functions directly for the properties that have inputs (or outputs) in the right-hand side. See \cite[Section IV.B]{BaM26} for a detailed discussion of this issue.

Next, assuming certain regularity of the flow, and using our converse UGAS Lyapunov results, we obtain a UGAS Lyapunov function w.r.t. this family of inputs, and then prove that it is an ISS Lyapunov function for the original system.

In our approach to converse ISS Lyapunov theory, we adopt the idea from \cite{SoW95}
of a reduction of the ISS analysis of the system to the UGAS analysis. 
However, while in the classical approach \cite{SoW95} one starts by transforming the original ISS system into an auxiliary UGAS system by closing the loop via a disturbed feedback (hence, changing both the system and in the input space), we just reduce the ISS property w.r.t. all admissible inputs to the question of UGAS of the original system with respect to the restricted class of inputs; see Lemma \ref{lem:ISStoROS}\ref{cond:ISStoROS4}. 
Hereby, we adapt the ideas from \cite{BaM26}, where we introduced the class of trajectory-dominated inputs to develop a converse Lyapunov theorem for boundedness of (finite-time) reachability sets \cite{MiW18b,BaM26}.
This allows to build a general converse Lyapunov theory for abstract systems given by the flow map even if we do not have equations defining the evolution of the system. 
However, to prove the regularity of the flow map of the system with respect to the family of admissibile input sets, we still employ the auxiliary feedback method, see Sections~\ref{sec:Converse ISS Lyapunov theorem for linear systems with unbounded input operators}, \ref{sec:Application to semilinear evolution equations and discussion}.

\emph{In Section~\ref{sec:Converse ISS Lyapunov theorem for linear systems with unbounded input operators}, we show that for general linear systems in Banach spaces with $\infty$-admissible input operators, ISS is equivalent to the existence of a Lipschitz continuous ISS Lyapunov function.} This fully solves the question about applicability of Lyapunov methods to the ISS analysis of linear boundary control systems: while coercive quadratic ISS Lyapunov functions exist for broad classes of systems with $2$-admissible input operators, this class is too restrictive to allow for $L^2$-ISS Lyapunov functions for systems with input operators that are not $2$-admissible. However, our results in this paper show that more general coercive ISS Lyapunov functions always exist.

Finally, in Section~\ref{sec:Application to semilinear evolution equations and discussion}, we apply our results to evolution equations with Lipschitz nonlinearities, and recover the main result from \cite{MiW17c}. 








\textbf{Notation:}
In the sequel, the nonnegative integers are denoted by $\N_0$, the natural numbers by $\N$, the real numbers by $\R$, the nonnegative real numbers by $\R_+$ and the open balls of radius $C$ around zero in normed vector spaces $X$, $U$ and $\U$, respectively, by $B_C$, $B_{C,U}$ and $B_{C,\U}$. For a subset $\Omega$ of a Banach space, we denote its closure by $\overline{\Omega}$.
We consider the standard classes of comparison functions (cf. \cite[p. xvi]{Mir23})
\begin{align*}
    \PP &\coloneqq \{\gamma: \R_+ \to \R_+ \,|\, \gamma(0) = 0,\ \gamma \text{ is continuous},  \\
    &\hspace{4cm}\gamma(s) > 0 \text{ for all } s > 0\}, \\
    \KK &\coloneqq \{\gamma \in \PP \,|\, \gamma \text{ is strictly increasing}\}, \\
    \KK_\infty &\coloneqq \{\gamma\in \KK \,|\, \gamma \text{ is unbounded}\}, \\
    \LL &\coloneqq \{\gamma: \R_+ \to \R_+ \,|\, \gamma \text{ is continuous and decreasing}  \\
    &\hspace{4cm}\text{with }\lim_{t \to \infty}\gamma(t) = 0\}, \\
    \KK\LL &\coloneqq \{\beta\colon \R_+ \times \R_+ \to \R_+ \,|\, \beta(\ph,t) \in \KK, \ \forall t \geq 0, \\
    &\hspace{4cm}\beta(r, \ph) \in \LL,\ \forall r > 0\}.
\end{align*}
We call a map $f\colon X \to Y$ between two Banach spaces $X,Y$ \emph{Lipschitz continuous on bounded sets}, if for each $C > 0$, there exists $L_C > 0$ such that for all $x_1,x_2 \in B_C$, it holds that
\begin{align*}
    \norm{f(x_1) - f(x_2)}_Y 
    &\leq L_C \norm{x_1 - x_2}_X.
\end{align*}


Moreover, for $k > 0$, let
\begin{align*}
    \operatorname{Lip}_k&\coloneqq \{\gamma\colon \R_+ \to \R_+ \,|\, \gamma \text{ is globally Lipschitz continuous} \\ &\hspace{2.7cm}\text{with Lipschitz constant at most } k\}.
\end{align*}

For $I \subset \R_+$ and a Banach space $Z$, we denote by $\LL^\infty(I, Z)$ the Lebesgue space of strongly measurable functions $f \colon I \to Z$ with a finite norm $\norm{f}_{\LL^\infty(I, Z)} \coloneqq \esssup_{t \in I}\norm{f(t)}_Z$ and 
by $\Cont(X,Y)$ the space of continuous maps between subsets of two normed vector spaces $X,Y$. The subspace of bounded functions in $\Cont(X,Y)$ is equipped with the supremum norm. By $\id$ we denote the identity map.

The following technical lemma from \cite[Lem. B.29]{Mir23} is useful.
\begin{lemma}\label{lem:boundSupremumDifference}
    Let $J$ be a nonempty index set and $v,w\colon J \to \R$. Then, the inequality
    \begin{align*}
        \sup\nolimits_{s \in J}\braces{v(s) - w(s)}
        &\geq\sup\nolimits_{s \in J}\!\braces{v(s)} - \sup\nolimits_{s \in J}\!\braces{w(s)}
    \end{align*}
    holds if the right-hand side is well-defined, i.e., if not both $\sup\nolimits_{s \in J}\!\braces{v(s)}$ and $\sup\nolimits_{s \in J}\!\braces{w(s)}$ are infinite.
\end{lemma}
Note that the case $\sup\nolimits_{s \in J}\!\braces{v(s)}= \infty > \sup\nolimits_{s \in J}\!\braces{w(s)}$ is not proven in \cite[Lem. B.29]{Mir23} but obvious.

\section{Abstract control systems}
\label{sec:Abstract control systems}

At first, we define the notion of an abstract control system.
\begin{definition}
\label{def:controlSystem}
    Consider a triple $\Sigma = (X, \U, \phi)$ consisting of
    \begin{enumerate}
        \item a normed vector space $(X, \norm{\ph}_X)$, called the \emph{state space}.
        \item a normed \emph{vector space of input values} $(U, \norm{\ph}_U)$ and the normed \emph{vector space of inputs} $\U = \LL^\infty(\I,U)$, equipped with the essential supremum norm, which we denote by $\norm \ph_\U$.
        \item\label{def:controlSystemProperty4} a map $\phi\colon D_\phi \to X$, $D_\phi \subset \I \times X \times \U$, called \emph{flow map}, so that for all $(x,u) \in  X \times \U$ there is a \emph{maximal existence time} $t_{\max} = t_{\max}(x,u) \in (0, + \infty]$ such that 
        \begin{center}
            $D_\phi \cap (\I \times \{(x,u)\}) = [0,t_{\max}) \times \{(x,u)\}$. 
        \end{center}
        The corresponding interval $[0,t_{\max})$ is called the \emph{maximal domain of definition} of the mapping $t \mapsto \phi(t, x, u)$, which we call a \emph{trajectory} of the system.
    \end{enumerate}
    
    The triple $\Sigma$ is called a \emph{(control) system} if it satisfies the following axioms.
    \begin{enumerate}
        \renewcommand\theenumi{($\Sigma$\arabic{enumi})}
        \item\label{cond:identityProperty} \emph{Identity property}: For all $(x,u) \in X \times \U$, it holds that $\phi(0, x, u) = x$.
        \item\label{cond:causalityProperty} \emph{Causality}: For all $(t,x,u) \in D_\phi$ and all $\widetilde u \in \U$ such that $u(s)= \widetilde u(s)$ for almost all $s \in [0,t]$, it holds that $[0,t] \times \{(x, \widetilde u)\} \subset D_\phi$ and $\phi(t,x,u) = \phi(t,x,\widetilde u)$.

        \item\label{cond:cocycleProperty} \emph{Cocycle property}: For all $x \in X$, $u \in \U$ and $t,s \geq 0$ 
        such that $t+s \in[0,t_{\max}(x,u))$, we have 
        \begin{align*}
            \phi(t + s, x, u) &= \phi\!\paren{s, \phi(t, x, u), u(\ph + t)}.
        \end{align*}
    \end{enumerate}
\end{definition}
The class of abstract control systems includes evolution equations with Lipschitz nonlinearities, sufficiently regular boundary control systems \cite[Sec. 4.1]{MiP20}, \cite{Sch20}, impulsive and switched systems, time delay systems, and broad classes of well-posed linear systems. 

In this setting, the admissible inputs are restricted to $\LL^\infty(\I, U)$, which is less general than the class of control systems that we considered in \cite{BDM26}. As discussed in \cite[Rem. 4.6]{JMP20}, if the inputs belong to $\LL^p$ spaces with $p<\infty$, as opposed to $\LL^\infty$, the definition of an ISS Lyapunov function has to be adapted.

For $u_1,u_2 \in \U$, denote their concatenation at time $t \in \I$ by
\begin{align}\label{eq:definitionConcatenation}
    u_1 \diamond_t u_2 (s)
    = \begin{cases}
        u_1(s) &\text{if } s \in [0,t),\\
        u_2(s-t) &\text{if } s \geq t.
    \end{cases}
\end{align}

\begin{remark}
Note that $u_1 \diamond_t u_2$ always belongs to $\U$.
Furthermore, $\Uc$ is invariant under time-shifts, and thus $\Uc$ satisfies the axiomatics as in \cite{MiW18b}.
\end{remark}

\begin{definition}
    \label{def:family-concatenation}    
Let for each $x \in X$ the \emph{set of admissible inputs} be given by $\U_x \subset \U$. We call $(\U_x)_{x\in X}$ a \emph{family of admissible input sets} if it satisfies the following properties:
    \begin{enumerate}
        \item $\mathcal U_x \neq \emptyset$ for every $x\in X$. 
    


        \item For any $x \in X$, $u\in \U_x$ and any $t \in [0,t_{\max}(x,u))$, we have that $u(t+\ph) \in \U_{\phi(t,x,u)}$.
    
        \item For any $x\in X$, we have
        \[
            u\in\mathcal U_x \land t \in [0,t_{\max}(x,u)) \land v\in\mathcal U_{\phi(t,x,u)} \implies u\diamond_t v\in\mathcal U_x.
        \]
    \end{enumerate}
\end{definition}

The cocycle property can be specialized to the family $(\U_x)_{x\in X}$ as follows:
\begin{proposition}
\label{prop:Alternative-cocycle}
Consider a control system $\Sigma$ and a family of admissible input sets
$(\U_x)_{x\in X}$. 
Then for all $x \in X$, $u \in \U_x$, $t \in [0,t_{\max}(x,u))$,  $v\in\U_{\phi(t,x,u)}$, and $s\ge0$ so that $t + s < t_{\max}(x,u\diamond_t v)$, we have 
    \begin{align}
    \label{eq:cocycle-alternative}
        \phi(t + s, x,  u \diamond_t v) &= \phi\!\paren{s, \phi(t, x, u), v}.
    \end{align}
\end{proposition}

\begin{proof}
    Item (iii) of Definition~\ref{def:family-concatenation} assures that $u \diamond_t v \in \U_x$. Applying the cocycle property with $u \diamond_t v$ instead of $u$ now gives the result.
\end{proof}

\begin{remark}
    Definition~\ref{def:family-concatenation} ensures that the system restricted to the family of admissible inputs is well-defined in the sense that all operations are compatible with the restriction to admissible inputs. In particular:
    \begin{enumerate}
        \item There is at least one admissible trajectory for each initial value.

        \item Applying the cocycle property on a trajectory with an admissible input leads to an initial condition for which the time-shifted input is admissible.
    
        \item The concatenation of admissible inputs along a trajectory is an admissible input. This property is crucial for the derivation of converse Lyapunov results.
    \end{enumerate}
\end{remark}

In what follows, a control system $\Sigma=(X,\U,\phi)$ endowed with a family of admissible input sets $(\U_x)_{x\in X}$ will be denoted by $(\Sigma, (\U_x)_{x\in X})$.

Control systems with a flow continuous with respect to time deserve special attention:
\begin{definition}\label{def:continuousControlSystem}
    We call a control system $(\Sigma, (\U_x)_{x\in X})$ a \emph{continuous control system}, if it satisfies
    \begin{enumerate}        \renewcommand\theenumi{($\Sigma$\arabic{enumi})}
        \setcounter{enumi}{3}
        \item\label{cond:continuity} \emph{Continuity}: For all $x\in X$ and all $u \in \U_x$, the map $t \mapsto \phi(t, x, u)$ is continuous on its maximal domain of definition.
    \end{enumerate}
\end{definition}

We need the following properties of control systems:
\begin{definition}\label{def:forwardComplete}
    We call a control system $(\Sigma, (\U_x)_{x\in X})$ \emph{forward complete (FC)}, if for each $x \in X$, $u \in \U_x$ and $t \in \I$ the value $\phi(t,x,u) \in X$ is well-defined.
\end{definition}

The next property is weaker than FC and ensures that a
trajectory can be prolonged as long as it stays bounded.

\begin{definition}[{\cite[Def.~1.4]{KaJ11a}}]\label{def:BIC}
    A control system $(\Sigma, (\U_x)_{x\in X})$ is said to satisfy the \emph{boundedness-implies-continuation (BIC) property} if for each $x \in X$, $u \in \U_x$ such that $t_{\max}(x,u)< \infty$, and for all $M > 0$, there exists  $t \in [0,t_{\max}(x,u))$ such that
$\norm{\phi(t,x,u)}_X \geq M$.
\end{definition}
The main stability notion for this work is
\begin{definition}\label{def:ISS}
    We call a control system $(\Sigma, (\U_x)_{x\in X})$ \emph{input-to-state stable (ISS)}, if it is forward complete and there exist $\beta \in \KK\LL$ and $\gamma \in \KK_\infty$ such that $\forall x \in X$, $\forall u \in \U_x$ the following holds:
    \begin{align}\label{ineq:ISS}
        \norm{\phi(t,x,u)}_X \leq \beta\!\para{\norm{x}_X,t} + \gamma\!\para{\norm{u}_\U}\!,\quad  t \in \R_+.
    \end{align}
\end{definition}

\begin{remark}
ISS of $(\Sigma, (\U_x)_{x \in X})$ coincides with
\begin{itemize}
    \item (classical) ISS, if for each $x \in X$, it holds that $\U_x = \U$;
    \item ISS with respect to small inputs 
    (see \cite[Definition 1]{CAI14}) 
    if there exists $R>0$ such that for every $x \in X$, it holds that $\U_x = \{u\in\Uc \,|\, \|u\|_{\Uc}\leq R\}$;
    \item uniform global asymptotic stability at zero (0-UGAS) if $\Uc_x=\{0\}$ for all $x\in X$.
\end{itemize}    
\end{remark}

We are now able to introduce ISS Lyapunov functions:
\begin{definition}
\label{def:ISSLyapunovFunction}
    Let $(\Sigma, (\U_x)_{x\in X})$ be a control system. We call $V \in \Cont\!\paren{X,\R_+}$ an \emph{ISS Lyapunov function 
    for $(\Sigma, (\U_x)_{x \in X})$} if it fulfills both of the following conditions:
	\begin{enumerate}
		\item\label{def:LyapFuncCondition1} There are $\alpha_1, \alpha_2 \in \KK_\infty$ such that for all $x \in X$
		\begin{align}
			\alpha_1\!\paren{\norm{x}_X} &\leq V(x) \leq \alpha_2\!\paren{\norm{x}_X}. \label{ineq:LFDefiniteness}
		\end{align}
        \item\label{def:LyapFuncCondition2} There exist $\chi \in \KK_\infty$ and $\psi \in \PP$ such that for all $x \in X$ and all $u \in \U_x$, it holds that
		\begin{align}			
            V\!\paren{x} \geq \chi\!\paren{\norm{u}_\U}
            \implies \dot V_u\!\paren{x} &\leq -\psi\!\paren{V\!\paren{x}},
            \label{ineq:LFFlow}
        \end{align}
        where we define the \emph{Lie derivative} $\dot V_u(x)$ by the \emph{upper right-hand Dini derivative} at zero for the function $t \mapsto V\!\paren{\phi(t,x,u)}$, i.e., 
        \begin{align*}
            \dot V_u(x)
            &\coloneqq \limsup_{s \searrow 0} \tfrac{1}{s}\paren{V\!\paren{\phi(s,x,u)} - V(x)}.
        \end{align*}
	\end{enumerate}
\end{definition}



In the seminal paper \cite{Son89}, Sontag has shown that the existence of an ISS Lyapunov function implies ISS for ODE systems. For infinite-dimensional systems, a counterpart of this result was shown in \cite[Theorem 1]{DaM13}. A generalization for general families of inputs is straightforward.
\begin{theorem}[Direct ISS Lyapunov theorem]
\label{thm:Lyapunov}
    Let $(\Sigma, (\U_x)_{x\in X})$ be a continuous control system satisfying the BIC property. If there exists an ISS Lyapunov function for $(\Sigma, (\U_x)_{x\in X})$, then $(\Sigma, (\U_x)_{x\in X})$ is ISS.
\end{theorem}

\section{Lyapunov characterization of global uniform asymptotic stability}
\label{sec:Discont-converseLyapunovThm}

Our main goal in this paper is to prove the Lipschitz continuous converse ISS Lyapunov theorem.  It is useful to start with a simpler concept, which can be called ISS with zero gain, but we prefer a more classic name:
\begin{definition}
\label{def:UGAS}
	We call a control system $(\Sigma, (\U_x)_{x\in X})$ \emph{uniformly globally asymptotically stable (UGAS)}, if it is forward complete and there exists $\beta \in \mathcal{KL}$ such that for all $x \in X$, $u \in \U_x$ and $t \in \R_+$, it holds that
	\begin{align}\label{ineq:UGAS}
		\norm{\phi(t,x,u)}_X \leq \beta\!\paren{\norm{x}_X\!, t}.
	\end{align}
\end{definition}

Our plan is first to derive converse Lyapunov results for UGAS, and then to reduce the ISS case to the UGAS case. 

\subsection{Discontinuous UGAS Lyapunov characterization}

Establishing regularity of a Lyapunov function is usually the main complexity in the derivation of converse Lyapunov results.
To simplify the problem, in this section we derive a \emph{discontinuous} exponential UGAS Lyapunov theorem without assuming any extra regularity from the flow $\phi$. These constructions will help us to derive Lipschitz continuous converse results later, under suitable regularity assumptions on $\phi$.

\begin{definition}
    \label{def:discontinuous-exp-UGAS-LF}
    A function $V:X \to \R_+$ is called a \emph{discontinuous exponential UGAS Lyapunov function}
    for $(\Sigma, (\U_x)_{x\in X})$ if the following conditions hold:
    
    \begin{enumerate}
        \item There exist $\alpha_1,\alpha_2\in\mathcal K_\infty$ such that
        \begin{equation}
        \label{eq:disc-LF-bounds}
            \alpha_1(\|x\|_X)
            \leq V(x)
            \leq \alpha_2(\|x\|_X),
            \qquad x\in X.
        \end{equation}
    
        \item There exists $a>0$ such that
        \begin{equation}
        \label{eq:disc-exp-decay}
         \hspace{-10mm}   V(\phi(t,x,u))
            \leq e^{-at}V(x),\quad x\in X,\  u\in \mathcal U_x, \  t\ge 0.
        \end{equation}
    \end{enumerate}
\end{definition}

\begin{remark}
The estimate (ii) easily implies that 
\[
\dot V_u(x) \leq - a V(x),\quad x\in X,\quad u\in \mathcal U_x.
\]
However, going back from this inequality to \eqref{eq:disc-exp-decay} requires certain regularity of the map $t\mapsto V(\phi(t,x,u))$. An easy sufficient condition would be the continuity of this map. 
\end{remark}






Now we can state a criterion for UGAS in terms of  discontinuous Lyapunov functions:
\begin{theorem}[UGAS Lyapunov characterization I]
    \label{thm:disc-exp-converse}    
    Let $(\Sigma, (\U_x)_{x\in X})$ be a forward complete control system.
    Then: 
    \begin{enumerate}
        \item\label{item:disc-exp-converse} If there is a \emph{discontinuous exponential UGAS Lyapunov function}
    for $(\Sigma, (\U_x)_{x\in X})$, then $(\Sigma, (\U_x)_{x\in X})$ is UGAS.
        \item[(ii)] 
        Let $(\Sigma, (\U_x)_{x\in X})$ be UGAS. Then, there is a \emph{discontinuous exponential UGAS Lyapunov function}
    for $(\Sigma, (\U_x)_{x\in X})$.
    
    More precisely, if there exist
        $\widetilde\alpha_1,\widetilde\alpha_2\in\mathcal K_\infty$ satisfying
        \begin{equation}
        \label{eq:Sontag-estimate-discontinuous}
            \widetilde\alpha_1
            \bigl(
                \|\phi(t,x,u)\|_X
            \bigr)
            \leq
            e^{-t}\widetilde\alpha_2(\|x\|_X)
        \end{equation}
        for every $x\in X$, $u\in\mathcal U_x$, and $t\geq0$, then, for
        every $a\in(0,1]$, a discontinuous exponential UGAS Lyapunov function is given by
        \begin{equation}
        \label{eq:W-discontinuous-converse}
            W_a(x)
            \coloneqq
            \sup_{v\in\mathcal U_x}
            \sup_{s\geq0}
            \left\{
                \widetilde\alpha_1
                \bigl(
                    \|\phi(s,x,v)\|_X
                \bigr)e^{as}
            \right\},\quad x \in X.
        \end{equation}
    \end{enumerate}
\end{theorem}

\begin{proof}
(i). Let $\alpha_1,\alpha_2\in\mathcal K_\infty$ and $a>0$ be as in
Definition~\ref{def:discontinuous-exp-UGAS-LF}. Then, for every
$x\in X$, every $u\in\mathcal U_x$, and every $t\geq0$,
\begin{align*}
    \alpha_1(\|\phi(t,x,u)\|_X)
    \leq
    V(\phi(t,x,u))
    &\leq
    e^{-at}V(x)
    \leq
    e^{-at}\alpha_2(\|x\|_X).
\end{align*}
Thus, $(\Sigma, (\U_x)_{x\in X})$ is UGAS with $\beta\in\mathcal{KL}$ given by
$  \beta(r,t)
    :=
    \alpha_1^{-1}
    \bigl(
        e^{-at}\alpha_2(r)
    \bigr)$, $r,t\ge 0$.


(ii). To prove the converse implication, assume that $(\Sigma, (\U_x)_{x\in X})$ is UGAS. Thus, there exists
$\beta\in\mathcal{KL}$ such that for every $x\in X$, $u\in\mathcal U_x$, and $t\geq0$
\[
    \|\phi(t,x,u)\|_X
    \leq
    \beta(\|x\|_X,t).
\]

By Sontag's $\mathcal{KL}$-lemma, there exist
$\widetilde\alpha_1,\widetilde\alpha_2\in\mathcal K_\infty$ such that
\eqref{eq:Sontag-estimate-discontinuous} holds. 
Fix $a\in(0,1]$ and define $W_a$ by
\eqref{eq:W-discontinuous-converse}.\\[-3mm]

\emph{Step 1: Sandwich inequality.} 
Pick any $x \in X$. 
Since $\mathcal U_x \neq \emptyset$, pick an arbitrary $v\in \U_x$. By the definition of $W_a$, and using the identity property of the flow, we obtain
\[
    W_a(x)
    \geq
    \widetilde\alpha_1(\|\phi(0,x,v)\|_X)
    =
    \widetilde\alpha_1(\|x\|_X).
\]

On the other hand, by
\eqref{eq:Sontag-estimate-discontinuous}, for every
$v\in\mathcal U_x$ and every $s\geq0$,
\begin{align*}
    \widetilde\alpha_1
    \bigl(
        \|\phi(s,x,v)\|_X
    \bigr)e^{as}
    &\leq
    e^{-s}\widetilde\alpha_2(\|x\|_X)e^{as}
    \\
    &=
    e^{-(1-a)s}\widetilde\alpha_2(\|x\|_X)
\leq
    \widetilde\alpha_2(\|x\|_X).
\end{align*}
Taking both suprema gives the sandwich inequality for $W_a$:
\[
    \widetilde\alpha_1(\|x\|_X)
    \leq
    W_a(x)
    \leq
    \widetilde\alpha_2(\|x\|_X),\quad x\in X.
\]
\vspace{-3mm}

\emph{Step 2: Exponential decay.} 
Fix $x\in X$, $u\in\mathcal U_x$, and $t\geq0$. By definition,
\begin{align*}
    W_a(\phi(t,x,u))
    &=
    \sup_{v\in\mathcal U_{\phi(t,x,u)}}
    \sup_{s\geq0}
    \left\{
        \widetilde\alpha_1
        \bigl(
            \|\phi(s,\phi(t,x,u),v)\|_X
        \bigr)e^{as}
    \right\}.
\end{align*}
By Definition~\ref{def:family-concatenation}, $u\diamond_t v\in\mathcal U_x$ for every
$v\in\mathcal U_{\phi(t,x,u)}$. 
Exploiting the cocycle property, we obtain
\begin{align*}
    &W_a(\phi(t,x,u)) =\\
    &e^{-at}
    \sup_{v\in\mathcal U_{\phi(t,x,u)}}
    \sup_{s\geq0}
    \left\{
        \widetilde\alpha_1
        \bigl(
            \|\phi(t+s,x,u\diamond_t v)\|_X
        \bigr)e^{a(t+s)}
    \right\}.
\end{align*}
As enlarging the sets over which the suprema are taken cannot make the supremum smaller, we obtain
\begin{align*}
    W_a(\phi(t,x,u))
    &\leq
    e^{-at}
    \sup_{w\in\mathcal U_x}
    \sup_{r\geq0}
    \left\{
        \widetilde\alpha_1
        \bigl(
            \|\phi(r,x,w)\|_X
        \bigr)e^{ar}
    \right\}
    \\
    &=
    e^{-at}W_a(x),
\end{align*}
which proves the claim.
\end{proof}

\begin{remark}
Our argument in Step 2 of the above proof may be viewed as a \emph{dynamic-programming argument}: 
evolving the system with a prescribed admissible input
during $[0,t]$ and optimizing only over the continuation after time
$t$ cannot produce a larger value than optimizing over all admissible
trajectories from the initial state.
\end{remark}

\subsection{Lipschitz continuous converse UGAS Lyapunov theorem}

Theorem~\ref{thm:disc-exp-converse} shows that for any UGAS control system, we can construct a Lyapunov function certifying UGAS without imposing any requirements on the regularity of the flow map $\phi$.



However, the continuity axiom \ref{cond:continuity} does not guarantee continuity of the constructed Lyapunov function $W_a$. 

Therefore, we adopt the notion of Lipschitz continuity on compact intervals with respect to trajectory-dominated inputs introduced in \cite{BaM26}.

\begin{definition}
\label{def:continuityOnCompactIntervalswrtInputs}
	Let $(\Sigma, (\U_x)_{x\in X})$ be a control system. We say that the \emph{flow of $(\Sigma, (\U_x)_{x\in X})$ is Lipschitz continuous on compact intervals with respect to $(\U_x)_{x\in X}$} if $(\Sigma, (\U_x)_{x\in X})$ is forward complete and 
    there exists a map $L: \R_+\times\R_+\to(0,+\infty)$ which is nondecreasing in both arguments such that for each $\tau > 0$ and each $C > 0$, for any $x_1, x_2 \in \ball_C$, and any $u_1 \in \U_{x_1}$, there exists $u_2 \in \U_{x_2}$, such that for all $t \in [0,\tau]$, it holds that
	\begin{align}\label{ineq:LipschitzContinuityOnCompactIntervalswrtInputs}
\hspace{-2mm}		\norm{\phi(t,x_1,u_1) - \phi(t,x_2,u_2)}_X \leq L(\tau,C)\norm{x_1 - x_2}_X\!.
	\end{align}
\end{definition}

\begin{remark}
\label{rem:regularityOfEvolutionEquationforRegularityOfSigma}
We refer to \cite{BaM26} for interpretations of Definition~\ref{def:continuityOnCompactIntervalswrtInputs} in terms of set-valued maps and for its relation to other concepts of Lipschitz continuity of trajectories.
\end{remark}

We introduce the following stronger notion of a  UGAS Lyapunov function.

\begin{definition}\label{def:UGAS-LF}
	Let $(\Sigma, (\U_x)_{x\in X})$ be a continuous control system.
    We call $V \in \Cont\!\paren{X,\R_+}$ a \emph{UGAS Lyapunov function} for $(\Sigma, (\U_x)_{x\in X})$, if it fulfills both of the 
    following conditions:
    \begin{enumerate}
		\item \label{cond:UGAS-LF1} There exist $\alpha_1, \alpha_2 \in \KK_\infty$ such that
		\begin{align}
			\alpha_1\!\paren{\norm{x}_X} &\leq V(x) \leq \alpha_2\!\paren{\norm{x}_X}, \quad x \in X.\label{ineq:UGAS-LFDefiniteness}
		\end{align}
		\item \label{cond:UGAS-LF2} There exists $\psi \in \KK_\infty$ such that for all $x \in X$ and $u \in \U_x$, it holds that
		\begin{align}
			\dot V_u\!\paren{x} \leq - \psi\!\paren{V\!\paren{x}}. \label{ineq:UGAS-LFFlow}
        \end{align}
	\end{enumerate}
\end{definition}

We will exploit the family of cut-off functions defined by
\begin{align*}
    G_k(z) \coloneqq\max\argbraces{0,z-\tfrac{1}{k}},\quad k \in \N, \quad z \in \R_+.
\end{align*}
Clearly, $G_k \in \operatorname{Lip}_1$ for every $k \in \N$:
\begin{align}
    \abs{G_k(z_1) - G_k(z_2)} \leq \abs{z_1 - z_2},\quad z_1,z_2\in\R_+. \label{ineq:estimateG_k1}
\end{align}

We proceed with a counterpart of the UGAS converse Lyapunov theorem \cite[Thm. 3.4]{KaJ11a}.
Our construction uses the Lyapunov functions $W_a$, which we have constructed earlier, as building blocks.

\begin{theorem}[UGAS Lyapunov characterization II]\label{prop:ROStoUGAS-LF}
    Let $(\Sigma, (\U_x)_{x\in X})$ be a continuous control system. Let 
    the flow of $\Sigma$ be Lipschitz continuous on compact intervals with respect to $(\U_x)_{x\in X}$. 
    
    
    Then, $(\Sigma, (\U_x)_{x\in X})$ is UGAS if and only if $(\Sigma, (\U_x)_{x\in X})$ admits a UGAS Lyapunov function  which is Lipschitz continuous on bounded sets.
\end{theorem}
\begin{proof}
    The proof of the direct Lyapunov theorem, i.e., $\exists$~UGAS Lyapunov function$\implies$UGAS, is completely analogous to the proof of Theorem \ref{thm:Lyapunov}.
    
    To prove the converse direction, we follow the approach in \cite[Thm. 3.4]{KaJ11a}.
    Let $(\Sigma, (\U_x)_{x \in X})$ be UGAS with $\beta$ as given in Definition \ref{def:UGAS}. Then, by Sontag's $\KK\LL$-Lemma \cite[Prop. A.7]{Mir23}, there exist $\widetilde\alpha_1, \widetilde\alpha_2 \in \KK_\infty$ such that for all $t\ge 0$, $x\in X$, $u\in \U_x$, we have
    \begin{align}
        \widetilde \alpha_1\!\paren{\norm{\phi(t,x,u)}_X} \leq e^{-t}\widetilde \alpha_2\!\paren{\norm x_X}. \label{ineq:estimateTildeAlpha1ofy-u}
    \end{align}

    We assume that $\widetilde \alpha_1 \in \operatorname{Lip}_1$ (otherwise take a smaller one with this property \cite[Lem. A.18]{Mir23}). 

    We define $\Theta = \Theta(R,q) \geq 0$ by
    \begin{align*}
        \Theta(R,q) \coloneqq \ln\!\paren{1 + q\widetilde \alpha_2(R)},\qquad R\geq 0,\ q\in\N.
    \end{align*}
    Now for all  $R\geq 0,\ q\in\N$, $x \in X$, $\norm{x}_X \leq R$, $u \in \U_x$ and $t \geq \Theta(R,q)$, the following inequality holds: 
    \begin{align}\label{ineq:UqBoundedInTime}
        \widetilde \alpha_1\!\paren{\norm{\phi(t,x,u)}_X}
        \leq e^{-t}\widetilde \alpha_2\!\paren{\norm x_X} \leq \tfrac 1 q.
    \end{align}
        
\noindent
    \emph{Step 1.1: Defining pre-Lyapunov functions:} We define for all $x \in X$ and $q \in \N$ the following maps
    \begin{align}
        U_q(x) 
        &\coloneqq \sup_{v \in \U_x} \sup_{s \geq 0}\!\big\{G_q\big(\widetilde\alpha_1\!\paren{\norm{\phi\argparen{s,x,v}}_X}\! \big)e^{s}\big\}. \label{eq:definitionUqUGAS}
    \end{align}
    The only difference between $U_q$ and $W_1$ is the map $G_q$, which makes the expression inside the supremum vanish whenever $\widetilde\alpha_1\!\paren{\norm{\phi\argparen{s,x,v}}_X} \leq \tfrac{1}{q}$.
    The reason for this adaptation is that this approach ensures Lipschitz continuity of the resulting Lyapunov function $V$ in a small neighborhood of the set $\braces{x \in X\,|\,V(x) = 0}$. Hence, several properties of $U_q$ can be verified as in Theorem~\ref{thm:disc-exp-converse}.



    \noindent
    \emph{Step 1.2: Sandwich bounds for $U_q$:} 
    This step is analogous to Step 1 in the proof of Theorem \ref{thm:disc-exp-converse} such that we obtain
    \begin{align*}
        G_q\argparen{\widetilde\alpha_1\!\paren{\norm{x}_X}} 
        &\leq U_q(x)\leq \widetilde \alpha_2\!\paren{\norm x_X},\quad x \in X,\quad q \in \N. 
    \end{align*}
\noindent
\emph{Step 1.3: Absolute bound for $U_q$:}
    Proceeding analogously to Step 2 of the proof of Theorem \ref{thm:disc-exp-converse}, we have for all $x \in X$, $u \in \U_x$ and $t \in \R_+$ that 
    %
    \begin{align}
        U_q(\phi(t,x,u)) 
        &\leq e^{-t} U_q(x). \label{ineq:boundU_q}
    \end{align}
    \noindent
    \emph{Step 1.4: Lipschitz continuity on bounded balls of $U_q$:}
    Pick any $q \in \N$, $R>0$ and let $x_1, x_2 \in B_R$. 
    By \eqref{ineq:UqBoundedInTime}, we have for $s \geq \Theta(\norm{x}_X\!,q)$ that
    \begin{align}
        G_q\big(\widetilde\alpha_1\!\paren{ \norm{\phi\argparen{s,x,u}}_X}\! \big) \leq \max\argbraces{0,\tfrac{1}{q}-\tfrac{1}{q}} = 0. \label{ineq:UqBoundedInTime2}
    \end{align}
    Hence, 
    %
    due to \eqref{ineq:UqBoundedInTime2}, we have 
    \begin{align}
\hspace{-2mm}        U_q(x) =& \sup_{v \in \U_x} \sup_{s \in [0,\Theta(R,q)]}\big\{G_q\big(\widetilde\alpha_1\!\paren{\norm{\phi\argparen{s,x,v}}_X}\! \big)e^{s}\big\}, \label{eq:definitionU_qBoundedDomain}
    \end{align}
    i.e., it is sufficient to evaluate $G_q\big(\widetilde\alpha_1\!\paren{\norm{\phi\argparen{s,x,v}}_X}\! \big)e^{s}$ on the compact time interval $[0,\Theta(\norm x_X\!,q)]$.
    
    Since the flow of $(\Sigma, (\U_x)_{x\in X})$ is Lipschitz continuous on compact intervals with respect to $(\U_x)_{x\in X}$, for these fixed $R,q$ there is a map $k_{x_1,x_2}\colon \U_{x_1} \to \U_{x_2}$ which maps each $u_1 \in \U_{x_1}$ to $u_2 = k_{x_1,x_2}(u_1) \in \U_{x_2}$ such that
    \begin{align*}
        \norm{\phi(t,x_1,u_1) - \phi(t,x_2,u_2)}_X \leq L(\Theta(R,q),\!R)\norm{x_1 - x_2}_X
    \end{align*}
    holds for all $t \leq \Theta(R,q)$ (if there are multiple such $u_2$ for given $x_1,x_2 \in X$, $u_1 \in \U_{x_1}$, the image $k_{x_1,x_2}(u_1)$ may be selected arbitrarily).

    Next we estimate $U_q(x_1) - U_q(x_2)$ using step-by-step the equality \eqref{eq:definitionU_qBoundedDomain}, the definition of $k_{x_1,x_2}$, Lemma \ref{lem:boundSupremumDifference} and \eqref{ineq:estimateG_k1}:
    \begin{align*}
    &U_q(x_1) - U_q(x_2)\\
        &=\sup_{u_1 \in \U_{x_1}} \sup_{s \in [0,\Theta(R,q)]} G_q \big(\widetilde\alpha_1\!\paren{\norm{\phi(s,x_1,u_1)}_X}\!\big)e^s  \\
        &\qquad\quad- \sup_{u_2 \in \U_{x_2}} \sup_{s \in [0,\Theta(R,q)]} G_q \big(\widetilde\alpha_1\!\paren{\norm{\phi(s,x_2,u_2)}_X}\!\big)e^s \nonumber\\
        &\qquad\leq \sup_{u_1 \in \U_{x_1}} \bigg\{\sup_{s \in [0,\Theta(R,q)]} G_q \big(\widetilde\alpha_1\!\paren{\norm{\phi(s,x_1,u_1)}_X}\!\big)e^s  \\
        &\qquad\quad- \sup_{s \in [0,\Theta(R,q)]} G_q \Big(\widetilde\alpha_1\!\paren{\norm{\phi(s,x_2,k_{x_1,x_2}(u_1))}_X}\!\Big) e^s\bigg\} \nonumber\\
        &\qquad\leq \sup_{u_1 \in \U_{x_1}} \sup_{s \in [0,\Theta(R,q)]} \Big|G_q \big(\widetilde\alpha_1\!\paren{\norm{\phi(s,x_1,u_1)}_X}\!\big)e^s \\
        &\qquad\qquad\qquad\qquad - G_q \Big(\widetilde\alpha_1\!\paren{\norm{\phi(s,x_2,k_{x_1,x_2}(u_1))}_X}\!\Big)e^s\Big| \nonumber \\
        &\qquad\leq \sup_{u_1 \in \U_{x_1}} \sup_{s \in [0,\Theta(R,q)]} \big|\widetilde \alpha_1\!\paren{\norm{\phi(s,x_1,u_1)}_X}  \\
        &\qquad\qquad \qquad\qquad-  \widetilde\alpha_1\!\paren{\norm{\phi(s,x_2,k_{x_1,x_2}(u_1))}_X}\!\big|e^s. \nonumber
    \end{align*}
    Employing that $\widetilde \alpha_1 \in \operatorname{Lip}_1$, the reverse triangle inequality, and the assumption that the flow of $(\Sigma, (\U_x)_{x\in X})$ is Lipschitz continuous on compact intervals with respect to $(\U_x)_{x\in X}$, we achieve
    \begin{align}
    U_q&(x_1) - U_q(x_2)\nonumber\\
        &\leq \sup_{u_1 \in \U_{x_1}} \sup_{s \in [0,\Theta(R,q)]} \big|\norm{\phi(s,x_1,u_1)}_X 
        - \norm{\phi(s,x_2,k_{x_1,x_2}(u_1))}_X\big| e^s\nonumber \\
        &\leq \sup_{u_1 \in \U_{x_1}} \sup_{s \in [0,\Theta(R,q)]} \!\!\!\!\!\! \!\norm{\phi(s,x_1,u_1) - \phi(s,x_2,k_{x_1,x_2}(u_1))}_Xe^s\nonumber \\
        &\leq \sup_{u_1 \in \U_{x_1}} \sup_{s \in [0,\Theta(R,q)]} L(\Theta(R,q),R) \norm{x_1 - x_2}_X e^s\nonumber \\
        &\leq L(\Theta(R,q),R)e^{\Theta(R,q)} \norm{x_1 - x_2}_X . \label{ineq:continuityUUOU_q(x)}
    \end{align}
    By exchanging the roles of $x_1$ and $x_2$, the analogous estimate holds for $U_q(x_2) - U_q(x_1)$. Thus, it follows for all $x_1,x_2 \in B_R$ that    
    \begin{align}
        |U_q(x_1) - U_q(x_2)|
        &= \max\argbraces{U_q(x_1) - U_q(x_2), U_q(x_2) - U_q(x_1)} \nonumber\\
        &\leq L(\Theta(R,q),R)e^{\Theta(R,q)} \norm{x_1 - x_2}_X, \label{ineq:propertyG2UUO} 
    \end{align}
    i.e., $U_q$ is Lipschitz continuous on bounded balls.

    \noindent
    \emph{Step 2.1: Construction of $V$:} Let the flow of $(\Sigma, (\U_x)_{x\in X})$ be Lipschitz continuous on compact intervals with respect to $(\U_x)_{x\in X}$ with a Lipschitz constant $L = L(\tau,C)$ 
    as in Definition~\ref{def:continuityOnCompactIntervalswrtInputs}.
    
    Then, by \eqref{ineq:propertyG2UUO}, the function $M\colon  (0,\infty) \times \N \to \R_+$, $M(R,q) \coloneqq L\!\paren{\Theta(R,q),R}e^{\Theta(R,q)}$ satisfies
    \begin{align}
        q \geq R &\implies M(q,q) \geq M(R,q) \label{ineq:propertyM1} \\
        x_1,x_2 \in B_R &\implies \abs{U_q(x_1) - U_q(x_2)} \leq M(R,q) \norm{x_1 - x_2}_X \label{ineq:propertyM2}
    \end{align}
    for all $q \in \N$ and $R > 0$.
    Then, we define
    \begin{align}
        V(x) \coloneqq \sum_{q = 1}^\infty \frac{2^{-q} U_q(x)}{1 + M(q,q)},\label{eq:LFHighRegularity}
    \end{align}
    which we show to be a UGAS Lyapunov function that is Lipschitz continuous on bounded sets:
     
    Indeed, one can prove analogously to \cite[Proof of Thm. 3.4, eq. (3.39) and following]{KaJ11a} and \cite[Proof of Thm. B.31, eq. (B.34) and following]{Mir23} the sandwich bounds and Lipschitz continuity on bounded balls of $V$. It remains to show the following:

        \noindent
    \emph{Step 2.2: Absolute bound for $V$:} By \eqref{ineq:boundU_q}, it holds that
    \begin{align}
        V\argparen{\phi(t,x,u)} \leq e^{-t}V(x).\label{ineq:exponentialDecreaseV}
    \end{align}
    The decay is then given by 
    \begin{align*}
        \dot V_u\!\paren{x}
        &= \limsup_{t \searrow 0}\tfrac{1}{t}\paren{V\!\paren{\phi(t,x,u)} - V\!\paren{x}} \\
        &\leq \limsup_{t \searrow 0}\tfrac{e^{-t} - 1}{t}V(x) 
        = - V\!\paren{x}.
    \end{align*}
    
    This completes the proof.
\end{proof}

\begin{remark}[Historical remarks]
Sontag's $\KL$-lemma has been proved in \cite[Proposition 7]{Son98}. It was immediately realized that it is ultimately helpful for the proof of converse Lyapunov results for UGAS \cite{TeP00, NTK99}. 
The use of cutoff maps $G_k$ for regularization of constructed Lyapunov functions goes back at least to Yoshizawa \cite[p. 98]{Yos59}.

For discrete-time systems, continuity of Lyapunov functions is not essential; hence, discontinuous converse Lyapunov results in this setting are quite classical, see \cite{NTK99, JiW02, GrK14}. 
In particular, Theorem~\ref{thm:disc-exp-converse} and the construction~\eqref{eq:W-discontinuous-converse} are continuous-time counterparts of converse Lyapunov results established for discrete-time systems in \cite[Theorem 6.2]{GrK14}. 
In fact, taking in the definition of a discontinuous UGAS Lyapunov function $t=1$ and noting that $e^{-a}<1$, we see that it can serve as an exponential UGAS Lyapunov function for discrete-time systems. With this interpretation, Theorem~\ref{thm:disc-exp-converse} recovers the corresponding result from \cite{GrK14}.

Theorem~\ref{prop:ROStoUGAS-LF}
is an extension of \cite[Thm. 3.4]{KaJ11a} to general sets of admissible inputs. The only significant difference is that we employ in our argument the discontinuous converse Lyapunov result, which simplifies the proof of 
Theorem~\ref{prop:ROStoUGAS-LF}. Otherwise, the proof is pretty close to \cite[Thm. 3.4]{KaJ11a}.
\end{remark}

\section{Converse ISS Lyapunov theorem}
\label{sec:Converse ISS Lyapunov theorem}

A classical way to prove converse Lyapunov results is to construct Lyapunov functions as suprema or integrals over a certain set of trajectories. But if the upper bound on this set of trajectories depends on inputs (as it does for ISS), then this approach becomes challenging, as Lyapunov functions are maps solely of a state (see \cite[Section IV.B]{BaM26} for a more detailed discussion).

To overcome this, Sontag and Wang in \cite{SoW95} have proposed a method to transform the original ISS system into an auxiliary system which is UGAS by inserting a specific disturbed feedback. As the right-hand side of a UGAS estimate does not depend on the input, this makes it possible to construct a UGAS Lyapunov function for this auxiliary system, and then show that it is an ISS Lyapunov function of the original system \cite{SoW95}. 
This method has been transferred to ISS of evolution equations in Banach spaces in \cite{MiW17c}.

This approach is, however, not applicable (at least in a straightforward way) to abstract systems defined via flow maps. Here, we follow the approach developed in \cite{BaM26} to derive a converse Lyapunov result for the \emph{bounded (finite-time) reachability sets (BRS)} property.
We avoid introducing auxiliary systems, and instead prove that ISS of $(\Sigma, (\U_x)_{x\in X})$ implies UGAS of $(\Sigma, (\U^\eta_x)_{x\in X})$ where $\U^\eta_x$ are the sets of trajectory dominated inputs.  
Then, we adapt the classic machinery including Sontag's $\KL$-Lemma to derive a Lipschitz continuous Lyapunov function for the system.

\subsection{Trajectory-dominated inputs}

We start by recalling the concept of trajectory-dominated inputs introduced in \cite{BaM26}.
\begin{definition}\label{def:TDI}
    Let $(\Sigma, (\U_x)_{x\in X})$ be a forward complete control system 
    and assume that $0 \in \U_x$ for each $x \in X$. 
    Let $\eta \in \KK_\infty$. For a given $x \in X$, we call the set
    \begin{align*}
        \U_x^\eta \coloneqq \Big\{v \in \U_x\,|
        \norm{v(\theta)}_U 
        \leq \eta\!\paren{\norm{\phi(\theta,x,v)}_X} \text{ for a.e. } \theta \ge 0\Big\},    
    \end{align*}
    the \emph{set of trajectory-dominated inputs with respect to initial value $x \in X$ and growth margin $\eta$}. We will use the colloquial term \emph{set of trajectory-dominated inputs} whenever $x$ and $\eta$ are clear from the context.
\end{definition}


The following lemma shows that the structure of the sets $(\U_x)_{x \in X}$ is inherited by the family $(\U^\eta_x)_{x \in X}$.
\begin{lemma}
\label{lem:concatenationTDI}
    Let $(\Sigma, (\U_x)_{x\in X})$ be a forward complete control system. Let $0 \in \U_x$ for each $x \in X$. Let $\eta \in \KK_\infty$ and $t \in \R_+$, $x \in X$ be given. The following holds: 
    \begin{enumerate}
        \item\label{cond:concatenationTDI1} $0 \in \U_x^\eta$, and thus $\U_x^\eta \neq\emptyset$.
        \item\label{cond:concatenationTDI2} For any $u\in \U_x^\eta$ and any $t\ge 0$, we have that $u(t+\ph) \in \U_{\phi(t,x,u)}^\eta$.
        \item\label{cond:concatenationTDI3} If $u \in \U_x^\eta$ and $v \in \U_{\phi(t,x,u)}^\eta$, then $u \diamond_t v \in \U_x^\eta$.
        \item\label{cond:concatenationTDI4} For any $u\in \U_x^\eta$ and $t \geq 0$, we have that $u \diamond_t 0 \in \U_x^\eta$.
    \end{enumerate}
    In particular, $(\U_x^\eta)_{x \in X}$ is also a family of admissible input sets 
    and $0 \in \U_x^\eta$ for each $x \in X$.
\end{lemma}
    
\begin{proof}
    \emph{\ref{cond:concatenationTDI1}:} Clear. 

    \emph{\ref{cond:concatenationTDI2}:} Since $\U_x^\eta \subset \U_x$, we have by Definition \ref{def:family-concatenation} that for each $t \geq 0$ and $u \in \U_x$ it holds that $u(t + \ph) \in \U_{\phi(t,x,u)}$.
    
    For any $u\in \U_x^\eta$ and any $t\ge 0$, we have 
    \[
    \norm{u(t+\theta)}_U \leq \eta\!\paren{\norm{\phi(t+\theta,x,u)}_X},\quad \text{for a.e. }\theta \ge 0.
    \]
    The cocycle property implies that
    \[
    \norm{u(t+\ph)(\theta)}_U \leq \eta\!\paren{\norm{\phi(\theta,\phi(t,x,u),u(t+\ph))}_X}
    \]
    for almost every $\theta \ge 0$. This shows that $u(t+\ph) \in \U_{\phi(t,x,u)}^\eta$.

    \emph{\ref{cond:concatenationTDI3}:}
    Let $t \in \R_+$, $x \in X$, $u \in \U_x^\eta$ and $v \in \U_{\phi(t,x,u)}^\eta$.
    Then, for a.e.  $\theta \geq 0$
    \begin{align*}
        \norm{(u \diamond_t v)(t + \theta)}_U 
        = \norm{v(\theta)}_U &\leq \eta\!\paren{\norm{\phi(\theta,\phi(t,x,u),v)}_X} \\
        &= \eta\!\paren{\norm{\phi(t + \theta,x,u \diamond_t v)}_X},
    \end{align*}
    where we used the cocycle property \ref{cond:cocycleProperty} in the last step. Due to causality \ref{cond:causalityProperty}, the bound from $u \in \U_x^\eta$ holds for almost all $\theta \in [0,t]$.
    Moreover, by Definition \ref{def:family-concatenation}, it holds that $u\diamond_t v\in\mathcal U_x$.
    Thus, we have
    $u \diamond_t v \in \U_x^\eta$.
    
    \emph{\ref{cond:concatenationTDI4}:} For any $u \in \U_x^\eta$ and $t \in \R_+$, the input $u \diamond_t 0$ is in $\U_x$ since $0 \in \U_{\phi(t,x,u)}$ by Statement \ref{cond:concatenationTDI1} and $(\U_x)_{x \in X}$ is closed under concatenations as shown in Statement \ref{cond:concatenationTDI3}. Obviously, the inequality
    \begin{align*}
        \norm{(u \diamond_t 0) (\theta)}_U 
        \leq \eta\!\paren{\norm{\phi(\theta,x,u \diamond_t 0)}_X}
    \end{align*}
    holds for a.e. $\theta \ge 0$.
\end{proof}

In the remainder, we will employ the following concept as an intermediate step to prove the existence of an ISS Lyapunov function for ISS systems.

\begin{definition}
\label{def:SRS}
    Let $(\Sigma, (\U_x)_{x\in X})$ be a forward complete control system 
    and $0 \in \U_x$ for each $x \in X$.
	We call $(\Sigma, (\U_x)_{x\in X})$ \emph{strongly robustly stable (SRS)}, if there exists a locally Lipschitz continuous margin $\eta \in \KK_\infty$ such that $(\Sigma,(\U_x^{\eta})_{x \in X})$ is UGAS.
\end{definition}

\subsection{Lyapunov criteria for ISS}
\label{sec:converseLyapunovThm}





With the notion of SRS at hand, we are now able to state the converse Lyapunov theorem.
    

\begin{theorem}[ISS Lyapunov characterization]
\label{thm:converseLyapunov}
	Let $(\Sigma, (\U_x)_{x\in X})$ be a continuous control system satisfying BIC 
    and $0 \in \U_x$ for each $x \in X$. 
    
    Consider the following statements:
    \begin{enumerate}
        \item\label{cond:ISSLyapunovFunction3} There exists an ISS Lyapunov function for $(\Sigma, (\U_x)_{x\in X})$ which is Lipschitz continuous on bounded sets.   
        \item\label{cond:ISSLyapunovFunction2} There is a continuous ISS Lyapunov function for $(\Sigma, (\U_x)_{x\in X})$. 
        \item\label{cond:ISSLyapunovFunction1} $(\Sigma, (\U_x)_{x\in X})$ is ISS.
        \item\label{cond:ISSLyapunovFunction4} $(\Sigma, (\U_x)_{x\in X})$ is SRS with some locally Lipschitz margin $\eta \in \KK_\infty$.
    \end{enumerate}
    Then the following implications hold:
    \begin{center}
    	\ref{cond:ISSLyapunovFunction3}$ \implies $\ref{cond:ISSLyapunovFunction2}$ \implies $\ref{cond:ISSLyapunovFunction1}$ \implies $\ref{cond:ISSLyapunovFunction4}.
    \end{center}
    Furthermore, if the flow of $(\Sigma, (\U_x^\eta)_{x\in X})$ is Lipschitz continuous on compact intervals with respect to $(\U_x^\eta)_{x\in X}$ with the same margin $\eta$ as in \ref{cond:ISSLyapunovFunction4}, then
    \begin{center}
    	\ref{cond:ISSLyapunovFunction4}$ \implies $\ref{cond:ISSLyapunovFunction3}.
    \end{center}
\end{theorem}

\begin{proof}
    We divide the proof into several partial results and proceed as depicted in Figure \ref{fig:convLyapunov}.

    \begin{figure}[htbp]
    \vspace{3mm}
    \noindent\hfill
    \begin{tikzpicture}[node distance =6mm and 14mm]
        \node (ISS) at (0,0) {ISS};
        \node[right=of ISS.east] (SRS) {SRS};
        \node[right=of SRS.east,align=center,text width=25.5mm] (UGAS-LF) {{$\exists$ Lipschitz\\ UGAS LF for $(\Sigma,(\U_x^\eta)_{x \in X})$}};
        \node[below=of UGAS-LF.south,align=center,text width=1.8cm] (LipISS-LF) {$\exists$ Lipschitz ISS LF};
        \node[align=center,text width=20mm] (contISS-LF) at (LipISS-LF-|ISS) {$\exists$ continuous ISS LF};
        \path
        (ISS) edge[thick,double,double equal sign distance,-{Implies[]}] node[anchor=south]{\footnotesize Lem.~\ref{lem:ISStoROS}} (SRS)
        (SRS) edge[thick,double,double equal sign distance,-{Implies[]}] node[anchor=south]{\footnotesize Thm. \ref{prop:ROStoUGAS-LF}} node[anchor=north]{\footnotesize $+$ regularity} (UGAS-LF.west)
        (LipISS-LF) edge[thick,double,double equal sign distance,-{Implies[]}] node[anchor=south]{\footnotesize clear} (contISS-LF)
        (contISS-LF) edge[thick,double,double equal sign distance,-{Implies[]}] node[anchor=east]{\footnotesize Thm.~\ref{thm:Lyapunov}} (ISS);
        \path (UGAS-LF) edge[every edge,thick,double,double equal sign distance,-{Implies[]}] node[anchor=west]{\footnotesize Lem. \ref{lem:UGAS-LFandConvLFtoISS-LF}} (LipISS-LF);

    \end{tikzpicture}%
    \caption{Outline of the proof of Theorem \ref{thm:converseLyapunov} for systems $(\Sigma, (\U_x)_{x \in X})$ which satisfy BIC.}
    \label{fig:convLyapunov}
\end{figure}
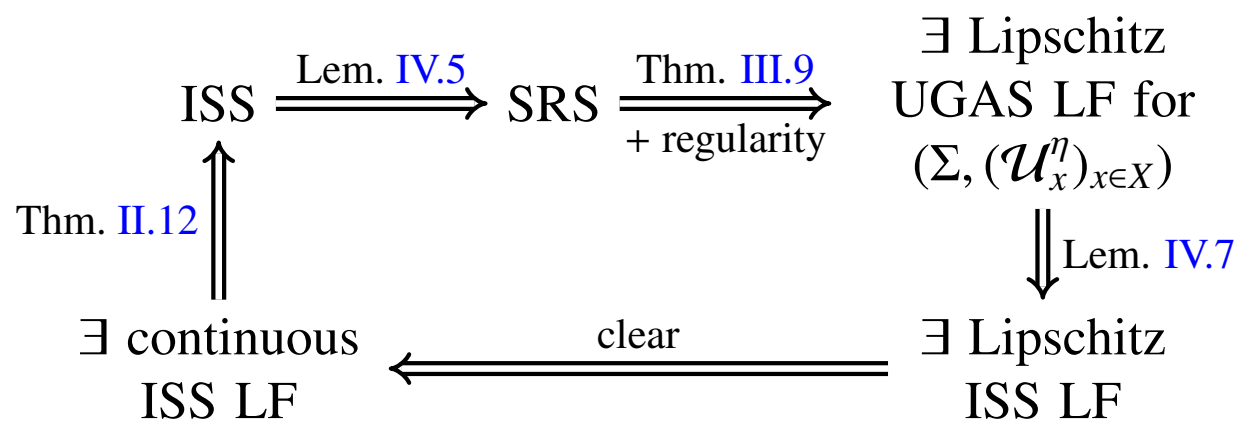%
    
    The implications \ref{cond:ISSLyapunovFunction3}$ \implies $\ref{cond:ISSLyapunovFunction2} and $ $\ref{cond:ISSLyapunovFunction2}$ \implies $\ref{cond:ISSLyapunovFunction1} follow from continuity and Theorem \ref{thm:Lyapunov}, respectively. Lemma~\ref{lem:ISStoROS} proves $ $\ref{cond:ISSLyapunovFunction1}$ \implies $\ref{cond:ISSLyapunovFunction4}. Under the additional regularity hypothesis, Theorem \ref{prop:ROStoUGAS-LF} and Lemma \ref{lem:UGAS-LFandConvLFtoISS-LF} prove \ref{cond:ISSLyapunovFunction4}$ \implies $\ref{cond:ISSLyapunovFunction3}.
\end{proof}

We start the proof of Theorem \ref{thm:converseLyapunov} by proving that ISS implies SRS:
\begin{lemma}
    \label{lem:ISStoROS}
    Let $(\Sigma, (\U_x)_{x\in X})$ be a forward complete control system 
    and $0 \in \U_x$ for each $x \in X$.
    Let $(\Sigma, (\U_x)_{x\in X})$ be 
    ISS with $\beta \in \KK\LL$, $\gamma \in \KK_\infty$ as in \eqref{ineq:ISS}. 
    Then for any locally Lipschitz continuous $\eta\in\Kinf$ such that    
    \begin{align}\label{ineq:upperBoundeta}
        \eta(s) \leq \gamma^{-1}\!\paren{\tfrac{1}{4}s},\quad s \in \R_+. 
    \end{align}    
    all of the following are true:
    \begin{enumerate}
        \item\label{cond:ISStoROS2} (Boundedness of $\U_x^\eta$).  There exists $\kappa \in \KK_\infty$ such that for each $x \in X$, it holds that 
        \[
        \U_x^\eta \subset 
        \braces{u \in \U_x \colon \norm{u}_\U \leq \kappa(\norm{x}_X)}.
        \] 
        
        \item\label{cond:ISStoROS4} $(\Sigma, (\U_x^\eta)_{x\in X})$ is UGAS.
        \item\label{cond:ISStoROS5} $(\Sigma, (\U_x)_{x\in X})$ is SRS with the above $\eta$.
    \end{enumerate}
    
\end{lemma}

\begin{remark}
    Note that Lemma \ref{lem:ISStoROS} \ref{cond:ISStoROS2} ensures that $\Sigma$ has bounded reachability sets for inputs in $\U_x^\eta$.
\end{remark}


\begin{proof}
    This proof is inspired by the approach in \cite[Lem. 3.3]{SoW99}. Let $\Sigma$ be ISS with $\beta \in \KK\LL$, $\gamma \in \KK_\infty$ as in \eqref{ineq:ISS}, and $\eta$ as in \eqref{ineq:upperBoundeta}. 
    Such $\eta$ is constructed in \cite[Lem. A. 18]{Mir23}. Let $\sigma \coloneqq \beta(\ph,0)$. Since $\beta(r,0) \geq r$ for every $r \geq 0$ by the identity property, we have $\sigma \in \KK_\infty$. For this $\eta$, consider $x \in X$ and $u \in \U_x^\eta$, i.e., such $u$ that
    \begin{align*}
        \norm{u(t)}_U \leq \eta\!\paren{\norm{\phi(t,x,u)}_X} \quad \text{for a.e. }t \in \R_+.
    \end{align*}

    
    \emph{Statement \ref{cond:ISStoROS2}:}    
    Let $T \in \R_+$ be arbitrary. Note that for any $u \in \U_x^\eta$, it holds that $u \diamond_T 0 \in \U_x^\eta$ by Lemma~\ref{lem:concatenationTDI}. By the ISS property and substituting $u \diamond_T 0$ instead of $u$, we obtain
    \begin{align*}
        &\norm{u}_{\LL^\infty([0,T),U)}  \leq \eta\!\paren{\beta\!\paren{\norm{x}_X,0} + \gamma\!\paren{\norm{u}_{\LL^\infty([0,T),U)}}}.
    \end{align*}
    
    From the last inequality and by definition of $\eta$, we see that if  $\gamma\!\paren{\norm{u}_{\LL^\infty([0,T),U)}} > \sigma\argparen{\norm{x}_X} = \beta\!\paren{\norm{x}_X,0}$, then, due to the strict monotonicity of~$\eta$, we have
    \begin{align*}
        \norm{u}_{\LL^\infty([0,T),U)} 
        &< \eta \!\paren{2\gamma\!\paren{\norm{u}_{\LL^\infty([0,T),U)}}} \\
        &\le \gamma^{-1}\!\paren{\tfrac{1}{4} \cdot 2\gamma\!\paren{\norm{u}_{\LL^\infty([0,T),U)}}} \\
        &\leq \norm{u}_{\LL^\infty([0,T),U)},
    \end{align*}
    and we obtain a contradiction. 
    Hence, we have $\gamma\!\paren{\norm{u}_{\LL^\infty([0,T),U)}} \leq \sigma\!\paren{\norm{x}_X}$, i.e., we arrive at the following bounds for $u$:
    \begin{align}
        \norm{u}_{\LL^\infty([0,T),U)}
        &\leq (\gamma^{-1} \circ \sigma) \argparen{\norm x_X}\!, \label{ineq:ISStoROSPreliminaryBound}
    \end{align}
    for arbitrary $T \in \R_+$. Hence, $\norm{u}_\U \leq \kappa(\norm{x}_X)$ for every $u \in \U_x^\eta$, where we set $\kappa = \gamma^{-1} \circ \sigma \in \KK_\infty$, just as claimed.

    \emph{Statement \ref{cond:ISStoROS4}:} 
    Since $\beta$ is a $\KK\LL$-function, for each $r > 0$, there exists $T_{r} \geq 0$, such that $\beta(r,t) \leq \frac{1}{4}r$ for all $t \geq T_{r}$. 
    We assume that $T = T_r$ is continuous and $T_r \geq 1$ for each $r > 0$. We show that such a choice is possible: By Sontag's $\KK\LL$-Lemma \cite[Prop. A.7]{Mir23}, there exist $\widetilde\alpha_1, \widetilde\alpha_2 \in \KK_\infty$ such that for all $r,t\ge 0$, we have
    \begin{align*}
        \widetilde \alpha_1\!\paren{\beta(r,t)} \leq e^{-t}\widetilde \alpha_2\!\paren{r}. 
    \end{align*}
    Obviously, we can choose $\widetilde \alpha_2$ such that $\widetilde \alpha_2(s) \geq e \widetilde \alpha_1(\tfrac{s}{4})$ for all $s \geq 0$. We define
    \begin{align*}
        T_r \coloneqq \ln \tfrac{\widetilde \alpha_2(r)}{\widetilde \alpha_1\argparen{\tfrac{r}{4}}}.
    \end{align*}
    Then, $T_r$ is continuous for $r \in (0,\infty)$ and satisfies $T_r \geq 1$.
    Moreover, Sontag's Lemma yields for every $r > 0$, $t \geq T_r$
    \begin{align*}
        \widetilde\alpha_1(\beta(r,t))
        \leq e^{-t}\widetilde\alpha_2(r)
        \leq \widetilde\alpha_1\argparen{\tfrac{r}{4}},
    \end{align*}
    i.e., $\beta(r,t) \leq \tfrac{r}{4}$ for every $r > 0$, $t \geq T_r$.
    
    If $x = 0$, Statement \ref{cond:ISStoROS2} gives $u = 0$, and ISS gives $\phi(t,0,0) = 0$. We may therefore fix $x \neq 0$
    and $u \in \U_x^\eta$. Set 
    $r_i \coloneqq 2^{1 - i}\sigma\!\paren{\norm{x}_X}$, $i \in \N_0$, $\Theta_0 = 0$, and for $i  \in \N$, recursively define $\Theta_i$ by 
    \begin{align*}
    \Theta_i \coloneqq\sum\nolimits_{j = 0}^{i - 1} T_{r_j} = T_{r_{i-1}} + \Theta_{i-1}.
    \end{align*}

    \noindent
    Claim: For each $i \in \N_0$, it holds that
    \begin{align}
        \sup\nolimits_{s \geq 0}\norm{\phi(\Theta_i + s,x,u)}_X
        &\leq r_i,  \label{ineq:ISStoROSoutput}
    \end{align}
    i.e., the trajectory after the time $\Theta_i$ stays in $\overline{B_{r_i}}$.
    
    \emph{Base case:} For $i = 0$, \eqref{ineq:ISStoROSoutput} follows from ISS of $\Sigma$ and \eqref{ineq:ISStoROSPreliminaryBound}:
    \begin{align*}
        \norm{\phi(t,x,u)}_X 
        &\leq \beta\!\paren{\norm{x}_X\!,t} + \gamma\!\para{\norm{u}_{\U}} \\
        &\leq 2\sigma\!\paren{\norm x_X} = r_0.
    \end{align*}
    
    \emph{Induction step:} Let \eqref{ineq:ISStoROSoutput} hold for $i \in \N_0$. 
    Using that $u \in \U_x^\eta$ and \eqref{ineq:ISStoROSoutput}, we obtain
    \begin{align}
        \norm {u(\ph + \Theta_i)}_\U 
        &= \esssup\nolimits_{s \geq 0}\norm{u(s + \Theta_i)}_U \nonumber\\
        &\leq \sup\nolimits_{s \geq 0} \eta\!\paren{\norm{\phi(s + \Theta_i,x,u)}_X} \nonumber\\
        &= \eta\!\paren{\sup\nolimits_{s \geq 0}\norm{\phi(s + \Theta_i,x,u)}_X} \nonumber\\
        &\leq \eta\!\paren{r_i}. \label{ineq:ISStoROSBoundInput}
    \end{align}
    
    As $\Theta_{i + 1} = \Theta_{i} + T_{r_i}$, the cocycle property \ref{cond:cocycleProperty} yields
    \begin{align*}
        \sup_{s \geq 0}&\norm{\phi(\Theta_{i + 1} + s,x,u)}_X \\
        &= \sup_{s \geq 0}\norm{\phi\!\paren{T_{r_i} + s,\phi\!\paren{\Theta_i, x,u},u(\ph + \Theta_i)}}_X.
    \end{align*}
    As $\Sigma$ is ISS, the estimates \eqref{ineq:ISStoROSBoundInput},  \eqref{ineq:ISStoROSoutput} imply 
    \begin{align*}
        \sup_{s \geq 0}&\norm{\phi(\Theta_{i + 1} + s,x,u)}_X \\
        &\leq \beta\!\paren{\norm{\phi\!\paren{\Theta_i, x,u}}_X\!, T_{r_i}} + \gamma \circ \eta\!\paren{r_i} \\
        &\leq \beta\!\paren{r_i, T_{r_i}} + \gamma \circ \eta\!\paren{r_i} \le \tfrac{1}{4} r_i + \tfrac{1}{4} r_i
        = \tfrac{1}{2} r_i 
        = r_{i + 1}.
    \end{align*}
    This proves the claim.



    From \eqref{ineq:ISStoROSoutput}, it follows for every $x \in X$ and $u \in \U_x^\eta$ that
    \begin{align*}
        \norm{\phi(t,x,u)}_X
        &\leq \widetilde \beta\!\paren{\norm x_X\!,t}\!,
    \end{align*}
    where $\widetilde \beta$ is chosen as an interpolation in time such as
    \begin{align*}
        \widetilde \beta\!\paren{r,t} 
        \coloneqq 2^{2 - i - \frac{t-\Theta_i}{\Theta_{i + 1} - \Theta_i}} \sigma\!\paren{r}
    \end{align*}
    for $t \in [\Theta_i, \Theta_{i + 1})$, $i \in \N_0$, i.e., system $\Sigma$ has the displayed decay bound. Here the times $\Theta_i = \Theta_i(r)$ are defined using $r_j = 2^{1 - j}\sigma(r)$ and we set $\widetilde \beta\!\paren{0,t} = 0$, $t \in \R_+$. This interpolation is continuous and decreasing in time but to ensure that it is strictly increasing in $r$ we introduce
    \begin{align*}
        \widehat \beta\!\paren{r,t} \coloneqq \max_{s \in [0,r]}\widetilde \beta\!\paren{s,t} + re^{-t}.
    \end{align*}
    Then it holds that $\widehat \beta \in \KK\LL$ and
    \begin{align*}
        \norm{\phi(t,x,u)}_X
        &\leq \widetilde \beta\!\paren{\norm x_X\!,t}
        \leq \widehat \beta\!\paren{\norm x_X\!,t}\!.
    \end{align*}
    Hence, $(\Sigma, (\U_x^\eta)_{x \in X})$ is UGAS. 
    
    \emph{Statement \ref{cond:ISStoROS5}}: directly follows from Statement \ref{cond:ISStoROS4}. 
\end{proof}

\begin{lemma}\label{lem:UGAS-LFandConvLFtoISS-LF}
    Let $(\Sigma, (\U_x)_{x\in X})$ be a forward complete continuous control system 
    and $0 \in \U_x$ for each $x \in X$.
    Let $\eta \in \KK_\infty$. Let the flow of $(\Sigma, (\U_x^\eta)_{x\in X})$ be Lipschitz continuous on compact intervals with respect to $(\U_x^\eta)_{x\in X}$. Let $(\Sigma, (\U_x^\eta)_{x\in X})$ be UGAS and possess a UGAS Lyapunov function $V$ which is Lipschitz continuous on bounded sets. 

    Then, $V$ is an ISS Lyapunov function for $(\Sigma, (\U_x)_{x\in X})$.
\end{lemma}




\begin{proof}
    Let $x \in X$ and $u \in \U_x$. 
    Let $V$ be a UGAS Lyapunov function which is Lipschitz continuous on bounded sets with $\alpha_1,\alpha_2,\psi \in \KK_\infty$ as in Definition \ref{def:UGAS-LF}.

    Obviously, \eqref{ineq:UGAS-LFDefiniteness} implies \eqref{ineq:LFDefiniteness}.

    Next, we define $\chi = \alpha_2\circ\eta^{-1} \circ ((1 + \eps)\id)$ for some $\eps > 0$. Fix $x \in X$ and $u \in \U_x$. If
    \begin{align*}
        V\!\paren{x} \geq \chi\!\paren{\norm{u}_\U}\!,
    \end{align*}
    then it holds by \eqref{ineq:LFDefiniteness} that
    \begin{align}\label{ineq:boundV1ofphi}
        \eta\argparen{\norm{x}_X} \geq (1 + \eps)(\norm{u}_\U).
    \end{align}
    
    \emph{Case 1 ($x \neq 0$):} Let us assume that $x \neq 0$, i.e., that $\norm{x}_X > 0$.
    Then, by \eqref{ineq:boundV1ofphi}, we have
    \begin{align*}
        \eta\argparen{\norm{x}_X}  > \norm{u}_\U.
    \end{align*}
     Due to continuity of $\phi(\ph, x, u)$, there exists some $\delta = \delta(x,u) > 0$ such that
    \begin{align*}
        \eta\argparen{\norm{\phi(t,x,u)}_X}  &> \norm{u}_\U &\forall t \in [0,\delta).
    \end{align*}
    We define $v \in \U_x$ by
    \begin{align*}
        v = u \diamond_\delta 0.
    \end{align*}
    By Definition \ref{def:family-concatenation} and since $0 \in \U_{\phi(\delta,x,u)}$, 
    it holds that $v \in \U_x^\eta$. Due to \eqref{ineq:UGAS-LFFlow}, $V$ satisfies
    \begin{align*}
        \dot V_v\!\paren{x} \leq - \psi\!\paren{V\!\paren{x}}
    \end{align*}
    and by the definition of the derivative and the causality property \ref{cond:causalityProperty}, it holds that
    \begin{align*}
        - \psi\!\paren{V\!\paren{x}} \geq \dot V_v\!\paren{x} 
        &= \limsup_{s \searrow 0} \tfrac{1}{s}\paren{V\!\paren{\phi(s,x,v)} - V(x)} \\
        &= \limsup_{s \searrow 0} \tfrac{1}{s}\paren{V\!\paren{\phi(s,x,u)} - V(x)} \\
        &= \dot V_u\!\paren{x}.
    \end{align*}
    Hence, for all $x \neq 0$, \eqref{ineq:LFFlow} is true. 
    
    \emph{Case 2 ($x = 0$):} It remains to show that
    \begin{align*}			
        V\!\paren{0} \geq \chi\!\paren{\norm{u}_\U}
        \implies \dot V_u\!\paren{0} &\leq -\psi\!\paren{V\!\paren{0}}.
    \end{align*}
    Indeed, the premise implies that $u \equiv 0$. Then, $u \in \U_x^\eta$ such that $\dot V_u\!\paren{0} \leq -\psi\!\paren{V\!\paren{0}}$ is a consequence of \eqref{ineq:UGAS-LFFlow}. 
    
    This concludes the proof of Lemma~\ref{lem:UGAS-LFandConvLFtoISS-LF}.
\end{proof}

\begin{remark}
    Based on the discontinuous converse UGAS Lyapunov theorem \ref{thm:disc-exp-converse}, we can prove a discontinuous converse Lyapunov theorem for ISS in the sense of \eqref{ineq:LFDefiniteness}--\eqref{ineq:LFFlow} using the method in this article. However, note that existence of a discontinuous function $V$ satisfying \eqref{ineq:LFDefiniteness}--\eqref{ineq:LFFlow} does not by itself
    give a direct ISS theorem because additional regularity along trajectories is required for trajectories with initial conditions and inputs satisfying $V(x) < \chi\argparen{\norm{u}_{\U}}$. In \cite{GrK14}, an analogous phenomenon was investigated for finite-dimensional discrete-time systems and appropriate results have been developed using more classical methods.
\end{remark}

\section{Converse ISS Lyapunov theorem for linear systems with unbounded input operators}
\label{sec:Converse ISS Lyapunov theorem for linear systems with unbounded input operators}

In this section, we prove that linear evolution equations with $\infty$-admissible input operators are ISS if and only if there exists a coercive ISS Lyapunov function which is Lipschitz continuous on bounded sets for such systems.


Consider a linear system
\begin{align*}
    \Sigma\colon \dot x(t) &= Ax(t) + Bu(t), \quad t \in \I,
\end{align*}
where $X, U$ are Banach spaces and the linear operator $A\colon D(A) \to X$, $D(A) \subset X$ generates a strongly continuous semigroup denoted by $(T(t))_{t \geq 0}$.

For a real $\lambda$ in the resolvent set of $A$, let $X_{-1}$ be the completion of $X$ in the norm $\norm{x}_{-1}=\norm{(\lambda I-A)^{-1}x}_X$, where $I$ is the identity operator. Write $T_{-1}(t)$ for the continuous extension of $T(t)$ to $X_{-1}$ and $L(U,X_{-1})$ for the space of bounded linear operators from $U$ to $X_{-1}$. 

We assume that 
 $B \in L(U,X_{-1})$ is an \emph{$\infty$-admissible} control operator for $(T(t))_{t \geq 0}$, i.e., 
\begin{align*}
t>0 \ \wedge \ u \in \LL^\infty(\R_+,U) \implies \int_0^t T_{-1}(t - s)Bu(s) \diff s \in X.
\end{align*}

As $B$ is $\infty$-admissible, a recent result \cite[Corollary 1.2]{APS26} shows that for any $x \in X$ and any $u\in\Uc$ the corresponding mild solution of $\Sigma$ 
\[
t\mapsto \phi(t,x,u)=T(t)x+\int_0^t T_{-1}(t - s)Bu(s) \diff s,\quad t\ge0,
\]
is a continuous function. 
With this, one can show that $\Sigma$ is a continuous forward complete control system in the sense of Definition~\ref{def:controlSystem}.

In this section, we consider the family of admissible input sets $\U_x = \U$ for every $x \in X$.

We call $B \in L(U,X_{-1})$ infinite-time $\infty$-admissible if there is $h>0$ such that for all $t>0$ and all $u \in \LL^\infty(\R_+,U)$ it holds that
    \begin{align*}
        \norm{\int_0^t T_{-1}(t - s)Bu(s) \diff s}_X \leq h \norm{u}_{\Uc}.
    \end{align*}
Note that if $(T(t))_{t \geq 0}$ is an exponentially stable semigroup, then $\infty$-admissibility and infinite-time $\infty$-admissibility are equivalent notions; see  \cite[Lem. 2.9]{JNPS18}.

The main result of this section is the Lyapunov characterization of ISS for general linear systems.
\begin{theorem}
\label{thm:Converse-Lyapunov-Theorem-linSys}
Let $B$ be $\infty$-admissible.
$\Sigma$ is ISS if and only if there is an ISS Lyapunov function which is Lipschitz continuous on bounded sets for $\Sigma$.
\end{theorem}

\begin{proof}
    "$\Leftarrow$". Follows from Theorem~\ref{thm:Lyapunov}.
    
    "$\Rightarrow$". We divide this proof into 3 parts.
    
    \noindent
    \emph{Step 1: $(\Sigma,(\U_x^\eta)_{x\in X})$ is UGAS for a suitable $\eta$:}
    
    Let $\Sigma$ be ISS. By \cite[Proposition 2.10]{JNPS18}, $T$ is exponentially stable, and $B$ is infinite-time $\infty$-admissible. Fix an admissibility constant $h > 0$ as above.
    Hence, there exist $M \geq 1$ and $a > 0$ so that for all $t\ge0$, $x\in X$ and $u\in\Uc$, we have
    \[
    \|\phi(t,x,u)\|_X \leq Me^{-at}\|x\|_X + h \norm{u}_{\Uc}.
    \]
    In particular, an ISS gain can be chosen as   $\gamma(r)=hr$, $r\geq 0$.
    
    
    
    
    Let $k = \frac{1}{4h}$ and $\eta \in \KK_\infty \cap \operatorname{Lip}_k$. The margin $\eta$ satisfies $\eta(r) \leq kr\leq\gamma^{-1}(r/4)$.
    Lemma~\ref{lem:ISStoROS} therefore ensures that $(\Sigma,(\U_x^\eta)_{x\in X})$ is UGAS.
    
    \noindent
    \emph{Step 2: Well-posedness of the auxiliary system:}
    
    Fix $d\in\D \coloneqq \overline{B_{1,\U}}$ and $x\in X$.
    On the Banach space of bounded continuous maps $\xi:\R_+\to X$, equipped with the supremum norm, define
    \begin{align}\label{eq:mildSolutionCL}
     \Phi_d(\xi)(t)&=T(t)x + \int_0^tT_{-1}(t-s)Bd(s)\eta\argparen{\norm{\xi(s)}_X} \diff s.
    \end{align}
    The input $s\mapsto d(s)\eta\argparen{\norm{\xi(s)}_X}$ is strongly measurable and bounded. By \cite[Corollary 1.2]{APS26} and the uniform admissibility estimate, $\Phi_d$ maps this Banach space into itself; exponential stability bounds its free term. Admissibility applied to the difference input, together with the $k$-Lipschitz property of $\eta$ and the reverse triangle inequality, gives
    \begin{align}
    \label{eq:Well-posedness-estimate}
        \|\Phi_d(\xi)&(t) - \Phi_d(\zeta)(t)\|_X \nonumber\\
        &= \Big\|T(t)x + \int_0^t T_{-1}(t - s) Bd(s)\eta\argparen{\norm{\xi(s)}_X} \diff s \nonumber\\
        &\qquad- T(t)x - \int_0^t T_{-1}(t - s) Bd(s)\eta\argparen{\norm{\zeta(s)}_X} \diff s\Big\|_X \nonumber\\
        &\leq \Big\|\int_0^t T_{-1}(t - s) Bd(s)(\eta\argparen{\norm{\xi(s)}_X} - \eta\argparen{\norm{\zeta(s)}_X})\diff s \Big\|_X \nonumber\\
        &\leq kh \esssup_{s\in[0,t]} \big(\|d(s)\|_U |\norm{\xi(s)}_X - \norm{\zeta(s)}_X |\big)
        \nonumber\\
        &\leq kh \esssup_{s\in[0,t]} |\norm{\xi(s)}_X - \norm{\zeta(s)}_X|
        \nonumber\\
        &\leq kh \sup_{s\in[0,t]} \norm{\xi(s) - \zeta(s)}_X.
    \end{align}
    Taking the supremum over $t \geq 0$ gives a contraction factor $kh\leq \frac{1}{4}$ in the supremum norm. The
    contraction theorem gives a unique bounded continuous fixed point, denoted by $\phi_{\CL}(\ph,x,d)$. It solves
    \begin{equation}\label{eq:CLlinearSystem}
     \Sigma_\CL:\quad\dot x(t)=Ax(t)+Bd(t)\eta\argparen{\norm{x(t)}_X}
    \end{equation}
    in the mild sense.
    
    \noindent
    \emph{Step 3: Lipschitz continuity of the flow with respect to $(\U_x^\eta)_{x \in X}$:} Let us show that $(\Sigma,(\U_x^\eta)_{x\in X})$ has a flow that is Lipschitz continuous on compact intervals with respect to $(\U_x^\eta)_{x\in X}$.
    
    Pick any $x_1,x_2 \in X$ and any $u_1 \in \U_{x_1}^\eta$ and define
    
    \begin{equation}\label{eq:implicitDefinitiond}
         d(t)=\begin{cases}
         u_1(t)/\eta(\norm{\phi(t,x_1,u_1)}_X),&\phi(t,x_1,u_1)\neq0,\\
         0,&\phi(t,x_1,u_1)=0.
         \end{cases}
    \end{equation}
    Then $d\in\D$. The bounded continuous trajectory $\phi(\ph,x_1,u_1)$ is a fixed point of \eqref{eq:mildSolutionCL} with $x:=x_1$, so it equals $\phi_{\CL}(\ph,x_1,d)$. Define, almost everywhere,
    \begin{equation}\label{eq:definitionOfu_2}
     u_2(t)=d(t)\eta(\norm{\phi_{\CL}(t,x_2,d)}_X).
    \end{equation}
    The fixed point is globally bounded, so $u_2\in\U$. Its mild equation gives $\phi(t,x_2,u_2)=\phi_{\CL}(t,x_2,d)$, and hence $u_2\in\U_{x_2}^\eta$. 
    
    Now using the estimate \eqref{eq:Well-posedness-estimate}, we proceed to 
    \begin{align*}
         &\norm{\phi(t,x_1,u_1) - \phi(t,x_2,u_2)}_X = \norm{\phi_\CL(t,x_1,d) - \phi_\CL(t,x_2,d)}_X \\
         &\leq \|T(t)\|\ \|x_1-x_2\|_X + kh \sup_{s\in[0,t]}\norm{\phi_\CL(s,x_1,d) - \phi_\CL(s,x_2,d)}_X\\
         &\leq  M\|x_1-x_2\|_X + kh \sup_{s\in[0,t]}\norm{\phi(s,x_1,u_1) - \phi(s,x_2,u_2)}_X.
    \end{align*}
    Taking the supremum over $[0,t]$ on both sides, we obtain 
    \[
    \sup_{s\in[0,t]}\norm{\phi(s,x_1,u_1) - \phi(s,x_2,u_2)}_X
    \leq \frac{M}{1-kh}\|x_1-x_2\|_X,
    \]
    which shows by employing Lemma~\ref{lem:ISStoROS} that $(\Sigma,(\U_x^\eta)_{x\in X})$ is UGAS and has a flow that is Lipschitz continuous on compact intervals with respect to $(\U_x^\eta)_{x \in X}$ with the same $\eta$. We remark that $\eta \in \KK_\infty \cap \operatorname{Lip}_k$ can be chosen arbitrarily.
    Thus Definition~\ref{def:continuityOnCompactIntervalswrtInputs} holds with $L(\tau,R) = \tfrac{M}{1 - kh}$, independent of $\tau$ and $R$.
    Application of Theorem~\ref{thm:converseLyapunov} concludes the proof.
\end{proof}

\section{Converse ISS Lyapunov theorem for semilinear evolution equations}
\label{sec:Application to semilinear evolution equations and discussion}

In this section, we specialize our results to \emph{semi-linear evolution equations} defined by
\begin{align*}
    \PDESigma\colon \dot x(t) &= Ax(t) + f\!\paren{x(t),u(t)}, \qquad t \in \I,
\end{align*}
where $X, U$ are Banach spaces, a linear operator $A\colon D(A) \to X$, $D(A) \subset X$ generates a strongly continuous semigroup over $X$ denoted by $(T(t))_{t \geq 0}$, and $f \colon X \times U \to X$ satisfies the following assumption: 

\begin{assumption}\label{as:regularityOff}
    The map $f \colon X \times U \to X$ is Lipschitz continuous on bounded balls, i.e., for each $C > 0$, there is $L_C > 0$ so that for all $x_1,x_2 \in B_C$, $u_1,u_2 \in B_{C,U}$, we have
    \begin{align*}
        \norm{f(x_1,u_1) - f(x_2,u_2)}_X\! \leq L_C \paren{\norm{x_1 - x_2}_X {+} \norm{u_1 - u_2}_U}.
    \end{align*}
\end{assumption}



Due to Lipschitz continuity on bounded balls of $f$, for each initial condition $x \in X$ and each input $u \in \U$, a unique maximal mild solution $\phi(\ph,x,u)$ of $\PDESigma$ exists \cite[Thm. 3.7]{Mir24}, i.e., there exists some $t_{\max} = t_{\max}(x,u) \in (0,  \infty]$ such that a function $\phi(\ph,x,u) \in \Cont([0,t_{\max}),X)$ satisfies for all $t \in [0,t_{\max})$ the integral equation
\begin{align}\label{eq:mildSolution}
    \phi(t,x,u) &=
    T(t)x + \int_0^t T(t - s) f\!\paren{\phi(s,x,u),u(s)} \diff s.
\end{align} 
ODEs are an important special case of the above systems, and thus our results are applicable in this context as well.


In this section, we consider the family of admissible input sets $\U_x = \U$ for every $x \in X$. In this case, for brevity, we write ``$\PDESigma$ is ISS'' instead of ``$(\PDESigma,(\U_x)_{x \in X})$ is ISS'' and proceed similarly for ISS Lyapunov functions and BIC. 


Converse ISS Lyapunov theorems for ODE systems and evolution equations $\PDESigma$ have been shown in  \cite{SoW95} and \cite{MiW17c}, respectively. In the following proposition we show how we can recover the main result from \cite{MiW17c} from our general framework.

\begin{proposition}
    Let $\PDESigma$ satisfy Assumption \ref{as:regularityOff}. Then, $\PDESigma$ is ISS if and only if  there is an ISS Lyapunov function for $\PDESigma$ which is Lipschitz continuous on bounded sets.
\end{proposition}
\begin{proof}
    We check the conditions of Theorem \ref{thm:converseLyapunov}.
    By \cite[Prop. 3.19]{Mir24},
    $\PDESigma$ is BIC and since we consider $\U_x = \U$ for all $x \in X$, we have $0 \in \U_x$ for all $x\in X$. Theorem \ref{thm:converseLyapunov} gives the direct implication. For the converse, let $\PDESigma$ be ISS. By Lemma~\ref{lem:ISStoROS}, $\PDESigma$ is SRS. Therefore, it remains to show that SRS of $\PDESigma$ implies that $\PDESigma$
    is Lipschitz continuous on compact intervals with respect to $(\U^\eta_x)_{x \in X}$, where $\eta \in \KK_\infty$ can be chosen to be the locally Lipschitz SRS margin.

    
    Hence, let $\eta \in \KK_\infty$ be such that $\PDESigma$ is SRS with margin $\eta$. We assume that $\eta \in \KK_\infty \cap \operatorname{Lip}_1$. Otherwise, we choose $\widetilde \eta \in \KK_\infty \cap \operatorname{Lip}_1$ such that $\widetilde \eta \leq \eta$ by \cite[Lem. A.18]{Mir23}. SRS of $\PDESigma$ with margin $\eta$ then implies SRS of $\PDESigma$ with margin $\widetilde \eta$ by definition.
    
    Next, we note that an evolution equation $\PDESigma$ which is SRS with margin $\eta$ is, in particular, robustly forward complete (RFC) with respect to trajectory-dominated inputs with growth margin $\eta$ as introduced in \cite[Def. II.8]{BaM26}. 
    
    Due to $\eta \in \KK_\infty \cap \operatorname{Lip}_1$, we have for each $C > 0$, for each $d \in D\coloneqq \braces{d \in U\,|\, \norm{d}_U \leq 1}$ and $x_1,x_2 \in B_C$ that $d \eta(\norm{x_1}_X), d \eta(\norm{x_2}_X) \in B_{C,U}$ due to global Lipschitz continuity with Lipschitz constant 1 of $\eta$, i.e. 
    \[
    \eta(s) = \abs{\eta(s) - \eta(0)} \leq \abs{s - 0} = s,\quad s \in \R_+.
    \]
    Assumption \ref{as:regularityOff}, $\eta \in \operatorname{Lip}_1$, and the reverse triangle inequality give
    \begin{align}
        &\|{f\argparen{x_1,d\eta\argparen{\norm{x_1}_X}} - f\argparen{x_2,d\eta\argparen{\norm{x_2}_X}}}\|_X \nonumber\\
        &\qquad\leq L_C\paren{\norm{x_1 - x_2}_X + \norm{d\eta\argparen{\norm{x_1}_X} - d\eta\argparen{\norm{x_2}_X}}_U} \nonumber\\
        &\qquad= L_C\paren{\norm{x_1 - x_2}_X + \norm{d}_U\abs{\eta\argparen{\norm{x_1}_X} - \eta\argparen{\norm{x_2}_X}}} \nonumber\\
        &\qquad\leq L_C\paren{\norm{x_1 - x_2}_X + \abs{\norm{x_1}_X - \norm{x_2}_X}} \nonumber\\
        &\qquad\leq 2L_C\norm{x_1 - x_2}_X. \label{ineq:regularityOff_feedback}
    \end{align}

    Moreover, by \eqref{ineq:regularityOff_feedback} and since $\PDESigma$ is RFC with respect to trajectory-dominated inputs with growth margin $\eta \in \KK_\infty \cap \operatorname{Lip}_1$,
    $(\PDESigma, (\U_x^\eta)_{x \in X})$ has a flow which is Lipschitz continuous on compact intervals with respect to $(\U_x^\eta)_{x\in X}$ with the same margin $\eta$ as was proven in \cite[Thm. V.1]{BaM26} (note that $\eta$ can be chosen the same for Statements (iii)--(v) of \cite[Thm. V.1]{BaM26}). 
    Therefore, SRS of $\PDESigma$ indeed implies that $\PDESigma$ is Lipschitz continuous on compact intervals with respect to $(\U^\eta_x)_{x \in X}$, where $\eta \in \KK_\infty$ can be chosen to be the locally Lipschitz SRS margin. 
    Application of Theorem~\ref{thm:converseLyapunov} finishes the proof.
\end{proof}

\begin{remark}
As ODEs with Lipschitz nonlinearities are a special case of evolution equations considered in this section, we obtain also Lipschitz continuous converse ISS Lyapunov theorem for ODE systems. However, in contrast to \cite{SoW95}, we cannot establish the existence of a smooth ISS Lyapunov function using our methods, as the smoothing technique used by \cite{SoW95} is not available in the infinite-dimensional setting.    
\end{remark}



\section{Discussion and conclusion}

We have proved a converse Lyapunov theorem for ISS of general infinite-dimensional control systems. 
As a corollary of our main result, we prove that linear well-posed ISS systems in Banach spaces always possess a Lipschitz continuous ISS Lyapunov function. 
This result will have broad applications in the control of boundary control systems, as many control techniques rely on the existence of coercive ISS Lyapunov functions, which our results ensure.
As another corollary, we recover the converse ISS Lyapunov theorem for evolution equations with Lipschitz nonlinearities proved in \cite{MiW17c}. 

As is usual for converse Lyapunov results, our construction is not explicit and relies on knowledge of the system flow, which is not given. 
Recently, in \cite{MiS26}, it was shown that for linear systems whose input operators are not 2-admissible, there are no quadratic coercive $L^2$-ISS Lyapunov functions. Possibly, one needs to consider Lyapunov functions generated by forms that are more general than bilinear forms. This is a great topic for future research. 

Another possible direction for future research is proving converse Lyapunov results for nonlinear boundary control systems. Our general result should apply to this case as well, provided one can show Lipschitz regularity of the flow map.

Finally, our new method for proving converse Lyapunov results for infinite-dimensional systems with inputs can be applied to many other stability and detectability properties, such as integral ISS, input-to-output stability, nonlinear robust detectability, etc. 







\bibliographystyle{elsarticle-num-names-nodoi}

\renewcommand{\bibfont}{\small}

\bibliography{Bibliothek,MyPublications}



\end{document}